\documentclass[]{elsarticle}
\pdfoutput=1

\makeatletter
\def\ps@pprintTitle{
\let\@oddhead\@empty
\let\@evenhead\@empty
\def\@oddfoot{\centerline{Preprint accepted to the Elsevier Journal of Computational and Applied Mathematics, Dec 2017.}}
\let\@evenfoot\@oddfoot}
\makeatother

\usepackage{hyperref}
\usepackage[margin=2.0cm]{geometry}
\usepackage{subcaption}
\usepackage{amsmath}
\usepackage{amssymb}
\usepackage{float}
\usepackage{multirow}
\usepackage{qtree}
\usepackage{algorithm}
\usepackage{color}
\usepackage[dvipsnames]{xcolor}

\begin{document}

\begin{frontmatter}
\title{Parallel accelerated cyclic reduction preconditioner for three-dimensional elliptic PDEs with variable coefficients}
\author[KAUSTaddress]{Gustavo Ch{\'a}vez\corref{correspondingAuthor}}
\cortext[correspondingAuthor]{Corresponding author}
\ead{gustavo.chavezchavez@kaust.edu.sa}

\author[AUBaddress]{George Turkiyyah}
\author[KAUSTaddress]{Stefano Zampini}
\author[KAUSTaddress]{David Keyes}
\address[KAUSTaddress]{King Abdullah University of Science and Technology, Thuwal, Saudi Arabia}
\address[AUBaddress]{American University of Beirut, Beirut, Lebanon}

\begin{abstract}
We present a robust and scalable preconditioner for the solution of large-scale linear systems that arise from the discretization of elliptic PDEs amenable to rank compression. The preconditioner is based on hierarchical low-rank approximations and the cyclic reduction method. The setup and application phases of the preconditioner achieve log-linear complexity in memory footprint and number of operations, and numerical experiments exhibit good weak and strong scalability at large processor counts in a distributed memory environment. Numerical experiments with linear systems that feature symmetry and nonsymmetry, definiteness and indefiniteness, constant and variable coefficients demonstrate the preconditioner applicability and robustness. Furthermore, it is possible to control the number of iterations via the accuracy threshold of the hierarchical matrix approximations and their arithmetic operations, and the tuning of the admissibility condition parameter. Together, these parameters allow for optimization of the memory requirements and performance of the preconditioner.
\end{abstract}


\end{frontmatter}

\section{Introduction}
This work focuses on the iterative solution of large-scale block tridiagonal linear systems of equations that arise from the discretization of elliptic partial differential equations on structured grids. Specifically, we demonstrate a parallel and scalable preconditioner based on an approximate factorization generated by the cyclic reduction algorithm~\cite{hockney65}. Cyclic reduction uses a sequence of Schur complement reduction steps, with each step eliminating half of the unknowns. While an exact cyclic reduction would result in prohibitively expensive dense matrix blocks, we exploit the data-sparsity of these resulting blocks by approximating them in a hierarchically low-rank form featuring log-linear storage.This work builds on \cite{Chavez2016ACR2D,Chavez2016ACR3D}, where a fast direct solver was introduced based on the synergy of parallel cyclic reduction and hierarchical matrices, and named accelerated cyclic reduction (ACR).

Iterative methods are advantageous for large-scale scientific computing since they feature tractable complexity and scalability, but their convergence is problem dependent. Direct methods, in contrast, guarantee a solution at the expense of higher complexity. Similar to the way in which incomplete factorizations such as the incomplete Cholesky factorization \cite{meijerink1977iterative} or the incomplete LU factorization \cite{saad2003iterative} accelerate the convergence, of Krylov methods, we propose a variable-accuracy ACR factorization that serves as a preconditioner to Krylov methods.

Since ACR is entirely algebraic, its range of applicability extends to problems with arbitrary coefficient structure, up to the amenability of rank compression. Furthermore, the ACR factorization and solve stages require only log-linear work and memory, which is particularly beneficial at large-scale. In addition, ACR exhibits substantial concurrency due to two separate characteristics: hierarchical matrix arithmetic operations expose substantial concurrency at the node level \cite{kriem05, Kriemann2014}, and the amount of distributed memory concurrency in cyclic reduction, which is based on red/black ordering of the block rows of the linear system, is proportional to the square root of the problem size in two-dimensions and to the cube root of the problem size in three-dimensions.

Numerical experiments document the robustness, performance, and memory consumption of the ACR preconditioner on a set of elliptic PDEs with heterogeneous coefficients. In particular, we study the variable-coefficient Poisson equation, the convection-diffusion equation, and the wave Helmholtz equation in heterogeneous media. For these equations, the numerical results show that ACR is a broadly applicable preconditioner that, without any equation-specific customizations. can be used to accelerate the convergence of Krylov methods very effectively even as the size of the linear system increases, the contract in the coefficients gets more pronounced, or the eigenvalue structure of the matris becomes more irregular. 

This rest of this paper is organized as follows. Section 2 reviews the literature on hierarchical matrices for the solution of elliptic PDEs. Section 3 reviews the basic elements of hierarchical matrix representations. Section 4 describes the accelerated cyclic reduction algorithm for the generation of the preconditioner, where hierarchical matrix representations are used for storing and manipulating the formally dense blocks that arise in the factorization process and shows the effect of tuning parameters on the hierarchical matrix structure. Section 5 describes the parallel version of the algorithm including the distributed memory parallelism of the overall factorization process  and the shared memory parallelism of the inner arithmetic operations on hierarchical matrices, and shows the weak and strong scalability of resulting algorithm. Sections 6 through 8 present detailed performance results on the effectiveness and near-optimal complexity of the ACR preconditioner on three problem categories of engineering interest: a variable coefficient Poisson problem, a non-symmetric convection-diffusion problem, and an indefinite Helmholtz equation. Section 9 presents the key conclusions.

\section{Literature review}
The last two decades have witnessed an increasing interest in the use of data-sparse approximations for the solution of linear systems. Leveraging an underlying hierarchically low-rank structure has been a successful strategy for improving the arithmetic complexity and memory footprint of direct solvers. As a result, direct solvers---as well as the closely related preconditioners obtained by aggressively truncating the rank of low-rank blocks---are becoming feasible candidates for tackling large-scale problems that traditional direct solvers are not able to handle due to memory requirements. In this section, we briefly discuss the two major directions towards fast direct solvers and preconditioners for the solution of sparse linear systems.

\subsection{Compression of dense frontal matrices}
The seminal work of Chandrasekaran et al.~\cite{chandrasekaran2010} showed that the off-diagonal blocks of the Schur complement of discretized elliptic PDEs can be efficiently represented with a hierarchical low-rank approximation. Using this property, methods such as the multifrontal solver~\cite{Duff83}, and other variants based on Schur complementation, can represent and perform arithmetic operations---of otherwise dense frontal matrices---using data-sparse formats.

An instance of the synergy of the multifrontal method with the hierarchical semiseparable (HSS) format \cite{Vandebril05} can be found in ~\cite{chandrasekarana2006, xia10,xia2010fast, xia2013randomized,Ghysels15,Wang2016}. However, other formats can be used to accelerate the multifrontal method, such as the hierarchical off-diagonal low rank (HODLR) format \cite{Ambikasaran2013} which lead to the multifrontal-HOLDR solver ~\cite{aminfar2016fast1}, or the block low-rank (BLR) format \cite{weisbecker2013improving} which lead to the multifrontal-BLR solver \cite{amestoy2015improving}, among others. For further discussion of the differences of each variant, we refer the reader to \cite{ChavezThesis2017}.

\subsection{Compression of the entire triangular factors}
An alternative technique that, rather than compressing individual blocks within the decomposition process, focuses on approximating the entire triangular factors as one hierarchical matrix has also been proposed. Instances of such strategy can be seen in the work of what is known as $\mathcal{H}$-Cholesky by Ibragimov et al.~\cite{ibra07} and $\mathcal{H}$-LU by Grasedyck et al.~\cite{blackBoxHLU08,grasedyck09}. The main idea is to create an $\mathcal{H}$-Matrix approximation of the sparse system with a clustering based on a nested dissection ordering of the unknowns. The nonzero blocks are approximated with a low-rank approximation, and an LU factorization is performed under the appropriate $\mathcal{H}$ arithmetic operations.

As with the previous section, different hierarchical formats can be used to approximate dense blocks. The work of Xia et al.~\cite{Xia10HSS_Cholesky} also proposes the construction of a rank-structured Cholesky factorization via the HSS hierarchical format, whereas the work of Pouransari et al.~\cite{Pouransari2016} approximates fill-in via low-rank approximations with the $\mathcal{H}^2$ format. We refer the reader to \cite{ChavezThesis2017} for a discussion of the differences of each strategy and a discussion of the implications of the choice of different hierarchical format regarding arithmetic operations count and memory requirements.

\section{Hierarchical matrix representations}
\label{sec:H}
A hierarchical matrix is a data-sparse representation that enables fast linear algebraic operations by using a hierarchy of off-diagonal blocks, each represented by a low-rank approximation or a small dense matrix, that can be tuned to guarantee a desired precision. The approximation, sometimes referred to as compression, is performed via singular value decomposition, or related methods that deliver low-rank approximations with fewer arithmetic operations than the traditional SVD method. For the representation to be effective in terms of arithmetic operations and memory requirements, the numerical rank must be significantly smaller than the sizes of the various matrix blocks that they replace.

\subsection{Overview of the $\mathcal{H}$-matrix format}

Formally, a hierarchical matrix in the $\mathcal{H}$-format \cite{hackbusch99,bebendorf2008_Book,hackbusch2015hierarchical}, can be constructed from four components: an index set, a cluster tree, a block cluster tree, and the specification of an admissibility condition.

\subsubsection{Index set}
The index set $\mathcal{I} = \{0,1,\dots,n-1\}$ represents the the nodal points of the grid under a certain ordering, such as the natural ordering.

\subsubsection{Cluster tree}
\begin{figure}[H]
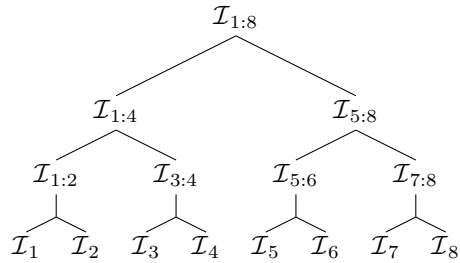

\Tree[.$\mathcal{I}$_{1:8} [.$\mathcal{I}$_{1:4} [.$\mathcal{I}$_{1:2} [ [.$\mathcal{I}$_{1} ][.$\mathcal{I}$_{2} ] ] ] [.$\mathcal{I}$_{3:4} [ [.$\mathcal{I}$_{3} ][.$\mathcal{I}$_{4} ] ] ] ]
[.$\mathcal{I}$_{5:8} [.$\mathcal{I}$_{5:6} [ [.$\mathcal{I}$_{5} ][.$\mathcal{I}$_{6} ] ] ] [.$\mathcal{I}$_{7:8} [ [.$\mathcal{I}$_{7} ][.$\mathcal{I}$_{8} ] ] ] ] ]
\caption{Binary cluster tree $\mathcal{T_I}$ of cardinality 8.}
\label{fig:tree}
\end{figure}

The cluster tree, denoted by $\mathcal{T_I}$, recursively subdivides the index set $\mathcal{I} \times \mathcal{I}$ until exhaustion. For simplicity, consider a binary cluster tree of cardinality 8 as shown in Figure \ref{fig:tree}.

\subsubsection{Block cluster tree}

\begin{figure}[H]
\centering
\includegraphics[width=0.75\textwidth]{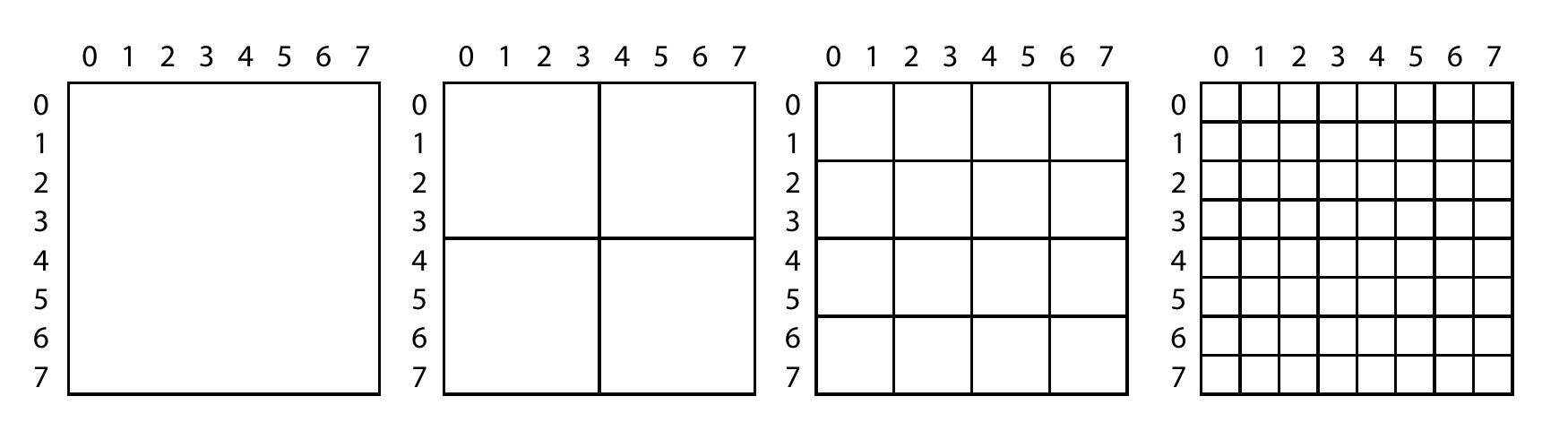}
\caption{Flat block partitioning at different levels $l$, without admissibility condition.}
\label{fig:Block}
\end{figure}

Once the cluster tree is defined, the block cluster tree maps matrix sub-blocks over the partitioning of the index set $\mathcal{I} \times \mathcal{I}$. An example of a clustering, frequently used by other data-sparse formats as we discuss in the next section, is a flat block-subdivision of the matrix in $l$ levels, as depicted in Figure \ref{fig:Block}. The $\mathcal{H}$-format, however, uses a discriminant to determine which blocks are further subdivided with the so-called admissibility condition.

\subsubsection{Admissibility condition}

\begin{figure}[H]
\centering
\includegraphics[width=0.75\textwidth]{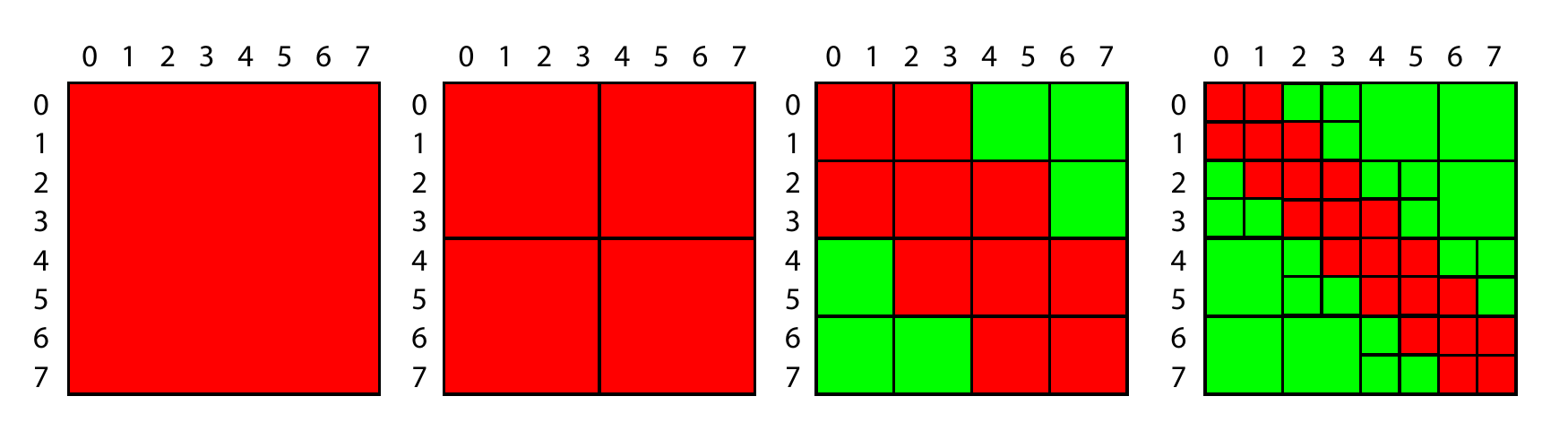}
\caption{Hierarchical block partitioning with standard admissibility condition. Green blocks are represented by low-rank approximations and red block with dense matrices.}
\label{fig:admis}
\end{figure}

Besides determining which blocks are further partitioned, the admissibility condition also determines which blocks are represented as a low-rank block (green) or a dense block (red), see Figure \ref{fig:admis}. 

A \emph{weak} admissibility criterion results in a coarse partitioning of the $\mathcal{H}$-matrix format where the off-diagonal blocks at every level of the hierarchy are all represented as low rank blocks. A \emph{standard} admissibility criterion allows a more refined blocking of the matrix given by the inequality:
\begin{equation}
min( diameter(\tau), diameter(\sigma) ) \leq \eta \cdot distance(\tau,\sigma)
\label{eq:admis}
\end{equation}

where $\tau$ and $\sigma$ denote two geometric regions defined as the convex hulls of two separate point sets $t$ and $s$ (nodes in the cluster tree). A matrix block $A_{ts}$ satisfying the previous inequality is represented in a low-rank form. The tuning parameter $\eta$ controls the weight of the distance function. Larger $\eta$ values admit larger blocks in the off-diagonal regions of the matrix as we will illustrate in Section 4.  

\subsubsection{Compression of low-rank blocks}

The last step for the construction of an $\mathcal{H}$-matrix is the choice of an algorithm to compute low-rank approximations for each of the blocks tagged as low-rank blocks as the product of two matrices of the from $UV^T$. Given a block of size $n \times n$, an effective compression leads to a tall and narrow matrix $U$ of size $n \times k$, and a short and wide matrix $V$ of size $k \times n$, where $k$ is the numerical rank of the block at some truncated accuracy $\mathcal{H}_{\epsilon}$. An effective compression means that the numerical rank $k$ is $k \ll n$. An efficient use of a hierarchical matrix to compress a given matrix has a balance between the numerical low-rank $k$ and a moderate number of low-rank blocks.

\subsection{Benefits of $\mathcal{H}$-matrix approximations}

$\mathcal{H}$-matrix approximations are especially useful for the particular class of matrices that arise from the discretization of elliptic operators with methods such as the boundary element method (BEM), finite-difference (FD), finite volumes (FV), or the finite element method (FEM). The resulting matrices and their Schur complements have blocks with bounded ranks that provide algorithmic gains while using $\mathcal{H}$-matrix storage and the set of algebraic operations that are available within the $\mathcal{H}$-format. In terms of storage, storing a dense matrix requires $\mathcal{O}(N^2)$ memory footprint, while its $\mathcal{H}$-matrix approximation counterpart can be stored in $O(N \log N)$ units of memory. For a comprehensive discussion of the construction of $\mathcal{H}$-matrices and their arithmetic operations, we refer the reader to \cite{hackbusch2015hierarchical}. 

\section{Accelerated cyclic reduction}
In this section we briefly review the cyclic reduction algorithm and describe the tunable accuracy accelerated cyclic reduction variant that improves its arithmetic and memory complexity estimates to near-optimal complexity for the variable-coefficient case.

\subsection{Cyclic reduction}

Cyclic reduction was introduced by Hockney in 1965 \cite{hockney65}, and then formalized by Buzbee and Golub in 1970 \cite{buzbee70}. Cyclic Reduction is a recursive algorithm for (block) tridiagonal linear systems. The algorithm consists of two phases: elimination and back-substitution. Elimination is equivalent to block Gaussian elimination without pivoting on a permuted system $(PAP^T)(Pu)=Pf$. The permutation matrix $P$ corresponds to a red-black ordering. 

\begin{figure}[ht!]
\centering
\includegraphics[width=0.75\linewidth]{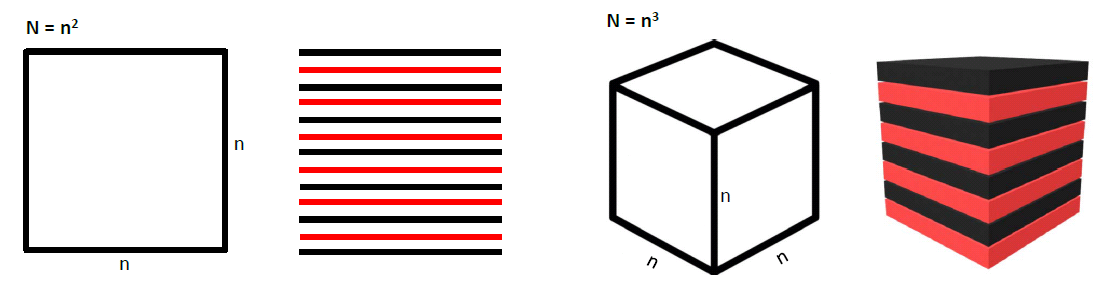}
\caption{
Plain grid in 2D and 3D, and its corresponding red/black coloring. In 2D $N=n^2$, each block row represents a line of size $n \times n.$, whereas in 3D $N=n^3$, each block row represents a plane of size $n^2 \times n^2.$}
\label{decomposition}
\end{figure}

The red-black ordering \textit{slices} the domain into lines or planes, depending on whether the underlying problem comes from a 2D or 3D problem respectively, as depicted on Figure \ref{decomposition}. This decomposition bears a similarity to the slice decomposition as reported in \cite{gugist2001}. This decomposition is also used in the sweeping preconditioner \cite{engquist2011sweeping,poulsonSweeping2013}.

\begin{figure}[ht!]
\centering
\includegraphics[width=\textwidth]{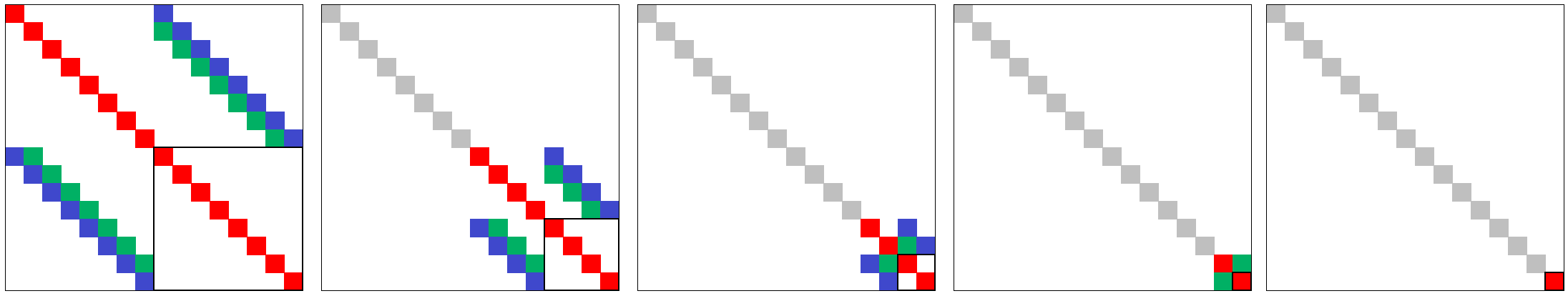}
\caption{Cyclic reduction preserves a block tridiagonal structure through elimination. Red blocks depict diagonal blocks, green blocks depict the innermost of the bidiagonal blocks, and blue blocks depict the outermost of the bidiagonal blocks. Gray blocks denote blocks in which elimination is completed.}
\label{fig:cr}
\end{figure}

Permutation decouples the system, and the computation of the Schur complement successively reduces the problem size by half. This process is recursive, and it finishes when a single block is reached, although the recursion can be stopped early if the system is small enough to be solved directly. The second phase performs a forward and backward substitution to find the solution. A graphical representation of the progression of elimination is shown in Figure \ref{fig:cr}. We refer to the reader to \cite{Chavez2016ACR2D} for an extended description of the cyclic reduction method.

\subsection{Accelerated cyclic reduction (ACR)}

In \cite{Chavez2016ACR3D}, we proposed the use of hierarchical matrices and their corresponding algebraic operations to improve on the computational complexity and the memory requirements of the classical cyclic reduction method, and named the resulting method accelerated cyclic reduction (ACR).

In generating the structure of the hierarchical matrix representations of the matrix blocks, we exploit the fact that, for a 3D problem, the domain is subdivided into $n$ planes each consisting of $n^2$ grid points. As a result, block rows of the matrix are identified with the planes of the discretization grid. We consider this geometry and use a two-dimensional planar bisection clustering when constructing each $\mathcal{H}$-matrix. In other words, ACR deals with $\mathcal{H}$-matrices with one dimension less than the original problem.

A standard admissibility condition was chosen, as opposed to a weak admissibility condition that the $\mathcal{H}$-matrix format also allows, because it provides the flexibility of selecting a range of coarser to finer blocks.

\subsection{Tuning parameters}
\label{sec:tuningETA}

There are three tuning parameters in the construction of an $\mathcal{H}$-matrix that can be leveraged to optimize memory requirements and performance: $\mathcal{H}_\epsilon$, $\eta$, and $n_{min}$. These, in turn, allow for a tunable accuracy ACR factorization which we use in this work as a preconditioner to Krylov methods. The cyclic reduction method was originally conceived as a direct solver; however, extensions of the use of CR as preconditioner have appeared in the literature \cite{rodrigue84, reusken00}, although to the best of our knowledge, none of them use hierarchical matrices.

The first parameter $\mathcal{H}_\epsilon$ controls the specified block-wise relative accuracy of the $\mathcal{H}$-matrix blocks tagged as low-rank. This parameter resembles the cut-off tolerance $\epsilon$ of the truncated SVD that disregards singular values to achieve an approximation accuracy of $\epsilon$.

\begin{figure}[H]
\begin{subfigure}{.5\textwidth}
\centering
\includegraphics[width=0.5\linewidth]{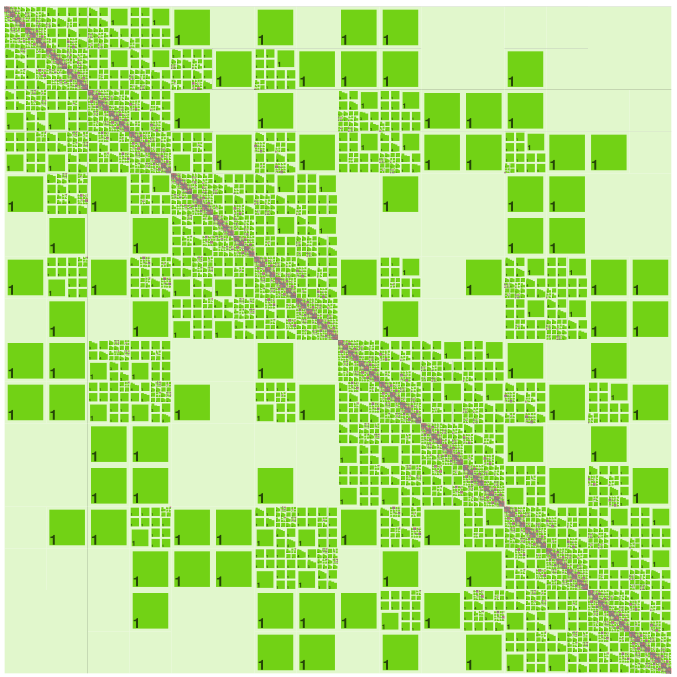}
\caption{$\eta$ = 2 (Strong admissibility)}
\label{fig:etaExperiment_a}
\end{subfigure}
\begin{subfigure}{.5\textwidth}
\centering
\includegraphics[width=0.5\linewidth]{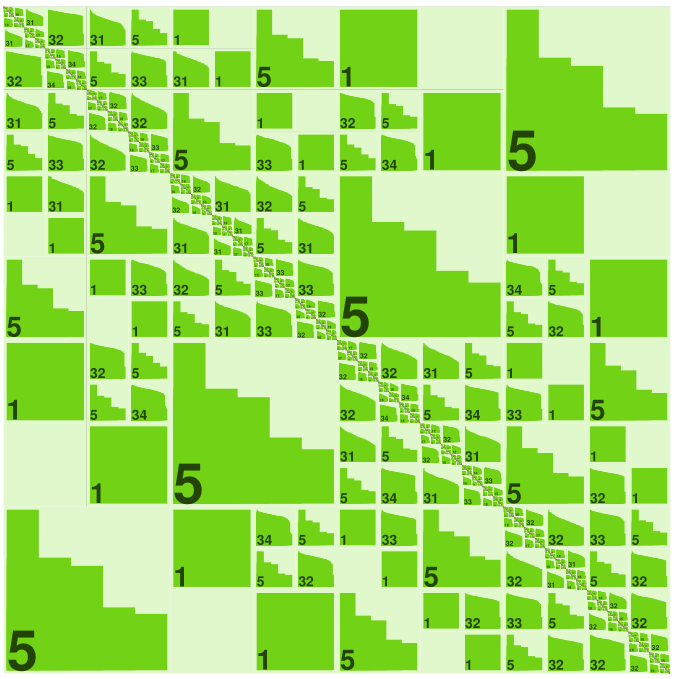}
\caption{$\eta$ = 64}
\label{fig:etaExperiment_b}
\end{subfigure}
\begin{subfigure}{.5\textwidth}
\centering
\includegraphics[width=0.5\linewidth]{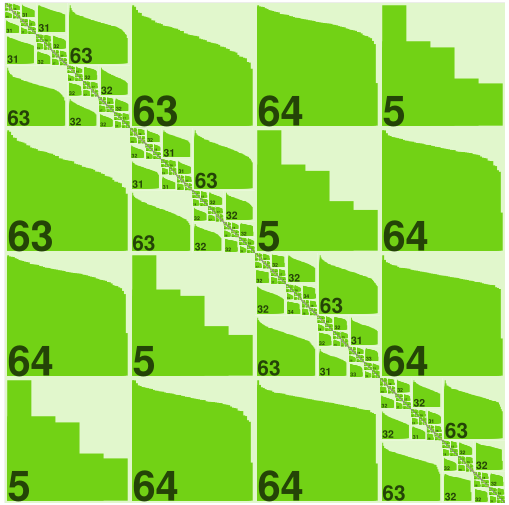}
\caption{$\eta$ = 128}
\label{fig:etaExperiment_c}
\end{subfigure}
\begin{subfigure}{.5\textwidth}
\centering
\includegraphics[width=0.5\linewidth]{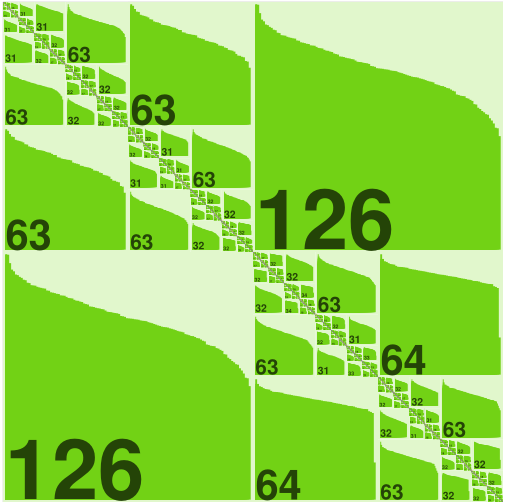}
\caption{$\eta$ = 256 (Weak admissibility)}
\label{fig:etaExperiment_d}
\end{subfigure}
\caption{$\mathcal{H}$-matrix structure for different parameter $\eta$ with fixed $\mathcal{H}_{\epsilon}$ and leaf size $n_{min}$=32. Matrix depicts a 2D variable coefficient Poisson problem with four orders of magnitude of contrast in the coefficient discretized with $N=128^2$ degrees of freedom. The numbers inside the green low-rank blocks denote the required numerical rank for the specified accuracy.}
\label{fig:etaExperiment}
\end{figure}

The second parameter is $\eta$, from the admissibility condition criterion. The case for choosing a standard admissibility condition (small $\eta$) is that, by further refining off-diagonals blocks, it is possible to achieve the same relative accuracy as with a weak admissibility (large $\eta$) but with smaller numerical ranks, albeit with more off-diagonal blocks (see Figure \ref{fig:etaExperiment_a} vs. \ref{fig:etaExperiment_d}). Numerical low-ranks are crucial to ensure economic memory consumption and overall high-performance.

Consider the computation of the approximate inverse in the $\mathcal{H}$-matrix format of a 2D variable-coefficient Poisson problem, with an error tolerance of three digits of accuracy in the Frobenius norm (${\left|\left|AA^{-1}-I\right|\right|_{F}}$), and a fixed accuracy parameter $\mathcal{H}_{\epsilon}$. The variable of interest in this experiment is the admissibility condition parameter $\eta$, which controls the block refinement as depicted in Figure \ref{fig:etaExperiment}.

Table \ref{table:etaTable} documents the memory requirements of each approximate inverse as a function of $\eta$. As shown, the optimal $\eta$ parameter resides in between strong admissibility ($\eta$=2) and weak admissibility ($\eta$=256). This tuning is a significant advantage for data-sparse formats that are not limited to the choice of weak admissibility, such as the $\mathcal{H}$-format. As the table shows, the most economic inverse regarding memory is not necessarily the representation with the smallest rank, since an aggressive refinement leads to a larger number of blocks and deeper cluster trees.

\begin{table}[H]
\centering
\begin{tabular}{|c|c|c|c|c|}
\hline
\textbf{$N$} & \textbf{$\eta$}  & \textbf{${\left|\left|AA^{-1}-I\right|\right|_{F}}$} & Max. Rank & Memory (Bytes) \\ \hline
$128^2$        & 2       & 5.0e-3     & 16      & 6.76e+7      \\ \hline
$\mathbf{128^2}$        & \textbf{64}      & \textbf{7.2e-3}     & \textbf{34}      & \textbf{6.64e+7}      \\ \hline
$128^2$        & 128     & 9.1e-3     & 64      & 8.64e+7      \\ \hline
$128^2$        & 256     & 9.1e-3     & 126     & 1.01e+8      \\ \hline
\end{tabular}
\caption{Memory consumption as a function of the tuning parameter $\eta$ for the computation of the approximate inverse in the $\mathcal{H}$-matrix format of a 2D variable-coefficient Poisson problem with four orders of magnitude of contrast in the coefficient, discretized with $N = 128^2$ degrees of freedom. Parameter $\eta$=2 depicts strong admissibility, while $\eta$=256 depicts weak admissibility; regarding memory requirements, $\eta$=64 is optimal.}
\label{table:etaTable}
\end{table}

\begin{figure}[H]
\centering
\includegraphics[width=0.4\linewidth]{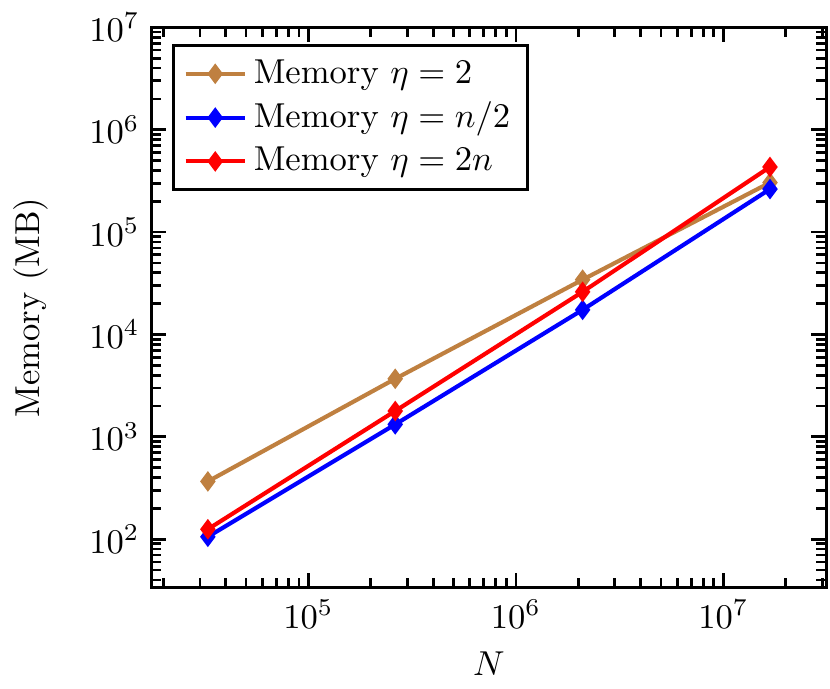}
\caption{Effect of tunable parameter $\eta$ on total memory consumption of ACR for a 3D variable-coefficient problem. The memory complexity estimate of $\mathcal{O}(N \log N)$ is achieved for $\eta=2$ which corresponds to strong admissibility. $\eta=2n$ which correspond to weak admissibility uses the most memory asymptotically, and an intermediate value of $\eta=n/2$ achieves the least amount of memory within the range of problem sizes considered.}
\label{fig:TuneEta}
\end{figure}

Since the memory consumption of ACR is determined by the sum of the memory consumption of each $\mathcal{H}$-matrix involved in elimination, an economical storage of each $\mathcal{H}$-matrix directly translates into savings to the overall ACR memory footprint. As shown in Figure \ref{fig:TuneEta}, tuning $\eta$ across a range of problem sizes has nuanced benefits. For linear systems in the order of a few millions of degrees of freedom a coarse block partitioning (close to weak admissibility) minimizes the overall memory consumption. However, for problems larger than a dozen of millions of unknowns block partitioning closer to strong admissibility is optimal to reduce memory requirements.

The third tuning parameter determines which blocks with less than or equal to $n_{min}$ rows or columns are stored as dense matrices, as it is more efficient to operate on them in dense rather than in low-rank form. It also alleviates unnecessarily deep binary trees in the structure of $\mathcal{H}$-matrices. Figure \ref{fig:etaExperiment} depicts these blocks in red.

\section{Hybrid distributed-shared parallelism}
The concurrency features of cyclic reduction have been evaluated in both distributed memory \cite{amodio2007algorithm, amodio1997cyclic, lin1994parallel, saghi1993predicting, sweet1988parallel}, and shared memory environments \cite{Strzodka2011cyclicGPU,quesada2011selecting,Zhang10}, although to the best of our knowledge, none of them use hierarchical matrices. We propose a hybrid model with MPI across the nodes and task-based parallelism across the cores in a node.

The distribution of parallel work was designed to accommodate the architecture of a modern supercomputer, such as the Shaheen Cray XC40 supercomputer at the King Abdullah University of Science $\&$ Technology. Shaheen is composed of 6,144 compute nodes, with each node holding 128GB of RAM and two Intel Haswell processors with 16 cores clocked at 2.3Ghz. The nodes are connected  with a Dragonfly network. All our reported numerical experiments were performed on this machine. 

Since the supercomputer architecture features multiple fast individual nodes, physically connected trough a high-speed network interconnect, we seek to maximize computation within nodes and minimize communication across nodes. Schur complementation is calculated locally within the nodes with a task-based parallel programming model. Dependencies to perform elimination and solve are fulfilled via a distributed memory programming model in which only the missing matrix blocks or vectors are provided with the message passing interface (MPI).

\subsection{Distributed memory parallelism}

\begin{figure}[ht!]
\centering
\includegraphics[width=0.5\linewidth]{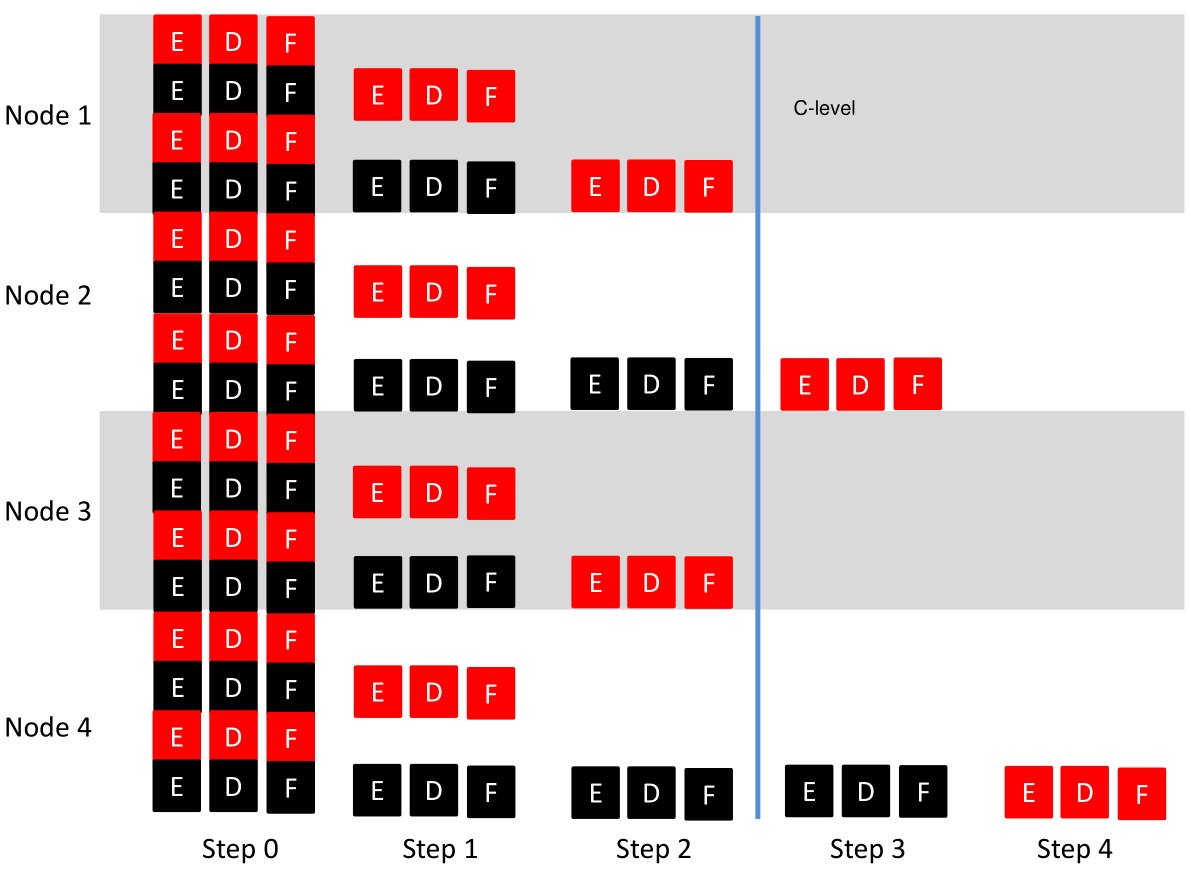}
\caption{Distribution of multiple planes per node for an example with $n=16$ planes and $p=4$ nodes.}
\label{fig:planesPerNode}
\end{figure}

Each plane of the computational domain is assigned to an MPI rank, but to minimize communication within the nodes, we allocate as many planes per node as memory allows. Let $p$ be the number of compute nodes and $n$ be the number of planes; therefore each node stores $n/p$ planes at the beginning of the factorization. Since ACR eliminates half of the planes at each step, after $r$ steps each node holds $n/(2^r p)$ planes. At level $r = \log (n/p)$, every node holds a single plane only. In the multigrid literature, this level is known as the C-level, i.e. the coarse level, illustrated in blue in Figure~\ref{fig:planesPerNode}. The remaining  $\log p$ steps beyond the C-level leave some compute nodes idle; fortunately, most of the remaining block operations have been completed by this step.

Distributed memory communication occurs only at inter-node boundaries, as the elimination of plane $i$ only requires planes $i-1$ and $i+1$. Thus up to the C-level, there are $O(p)$ messages per step, each transmitting planes of size $O(k~n^2 \log n)$. Beyond the C-level, there are $O(p/2+\cdots+1) \approx O(p)$ communication messages, adding up to a total communication volume of $O(k~p~n^2 \log n~ ( \log \frac{n}{p}+1 ) )$, where $n$ is the size of the linear dimension and $k$ is the rank of the low-rank approximation.

\subsection{Shared memory parallelism}

We exploit the concurrency of hierarchical matrix algebra at node level \cite{kriem05}, in particular, through a task-based programming model \cite{Kriemann2014}. The HLibPro package relies on the Intel Threading Building Blocks library \cite{pheatt2008intel} to build directed acyclic graphs for the dependencies between tasks, which might involve recursion. The allocation of cores to perform hierarchical matrix algebra within the node depends on the number of planes per node.

For instance, for a node with thirty-two processors and four planes per node, we set four MPI processes, and for each plane we allocate eight processors to perform task-based parallelism. Resource allocation for either block row processing (communication of planes) or parallel task-based hierarchical arithmetics (computation of Schur complement) can be tuned to maximize performance or memory availability.

We refer to the reader to \cite{Chavez2016ACR3D} for an extended discussion of the parallel features of ACR, and the derivation of its parallel complexity estimates.

\subsection{Parallel scalability}

The parallel scalability of the ACR preconditioner is demonstrated below on a constant-coefficient Poisson equation with homogeneous Dirichlet boundary conditions in the unit cube, discretized with the 7-point finite-difference star stencil. We use the conjugate gradient method \cite{hestenes1952methods} with a convergence criterion in the 2-norm of the relative residual down to $10^{-8}$. This problem results in a symmetric positive definite matrix whose factors exhibit rapid decay of the singular values of off-diagonal blocks and offers an ideal testbed to show the parallel scalability of the hierarchically low-rank algorithmic computations. The following experiments improve on previous timings for the setup phase performed at smaller accuracies \cite{Chavez2016ACR3D}, demonstrating that the use of ACR as a preconditioner with looser tolerances, i.e. smaller ranks, also benefits scalability.

\subsubsection{Strong scaling}

\begin{figure}[H]
\centering
\begin{subfigure}{0.4\textwidth}
\centering
\includegraphics[width=0.85\linewidth]{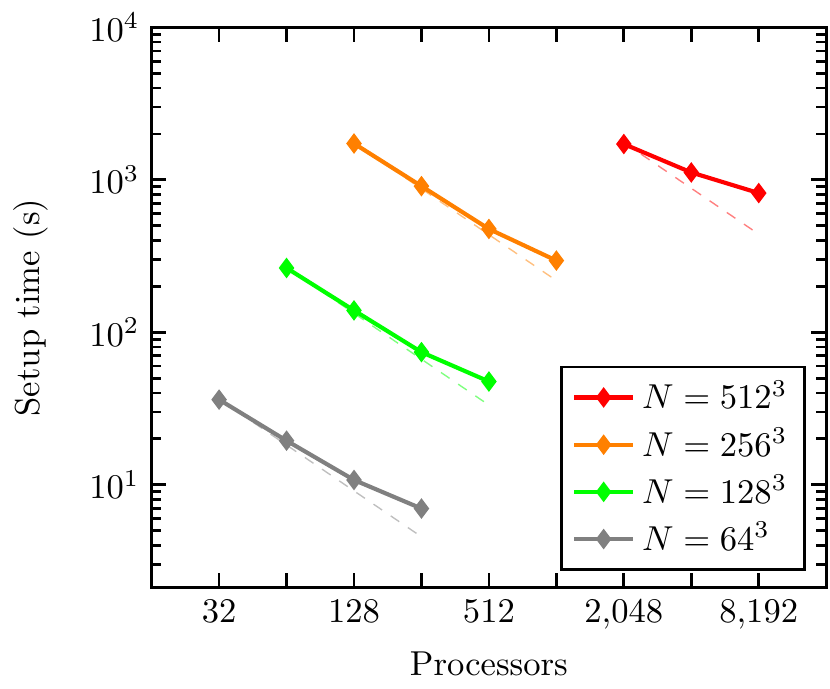}
\caption{Strong scaling of the preconditioner setup.}
\label{fig:strong1}
\end{subfigure}
\begin{subfigure}{0.4\textwidth}
\centering
\includegraphics[width=0.85\linewidth]{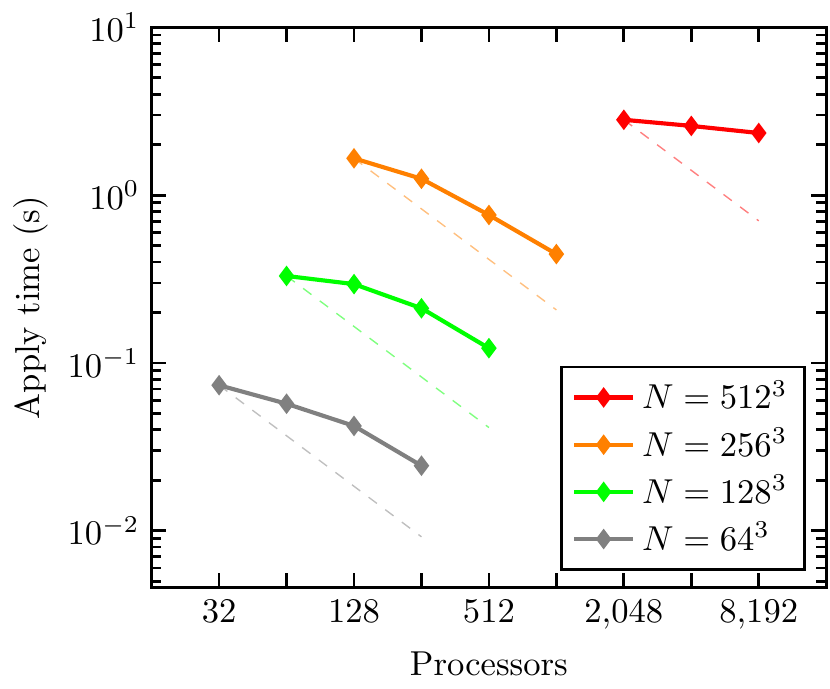}
\caption{Strong scaling of the preconditioner apply.}
\label{fig:strong2}
\end{subfigure}
\caption{ACR preconditioner strong scalability for the solution of the constant-coefficient Poisson equation.}
\end{figure}

Figures \ref{fig:strong1} and \ref{fig:strong2} show the total time in seconds for the setup and solve per iteration of the ACR preconditioner in a strong scaling setting; dashed lines indicate ideal scaling, a reduction in time by a factor of two as we double the number of processors.

The most time-consuming phase of the ACR preconditioner benefits the most as the number of processors increases for a variety of problem sizes. Nonetheless, the ideal scaling of the solve stage deteriorates at large processor counts as factors such as hardware latency play a significant role in this computationally lightweight kernel solely based on $\mathcal{H}$ matrix-vector multiplications.

\subsubsection{Weak scaling}

\begin{figure}[H]
\centering
\begin{subfigure}{0.4\textwidth}
\centering
\includegraphics[width=0.85\linewidth]{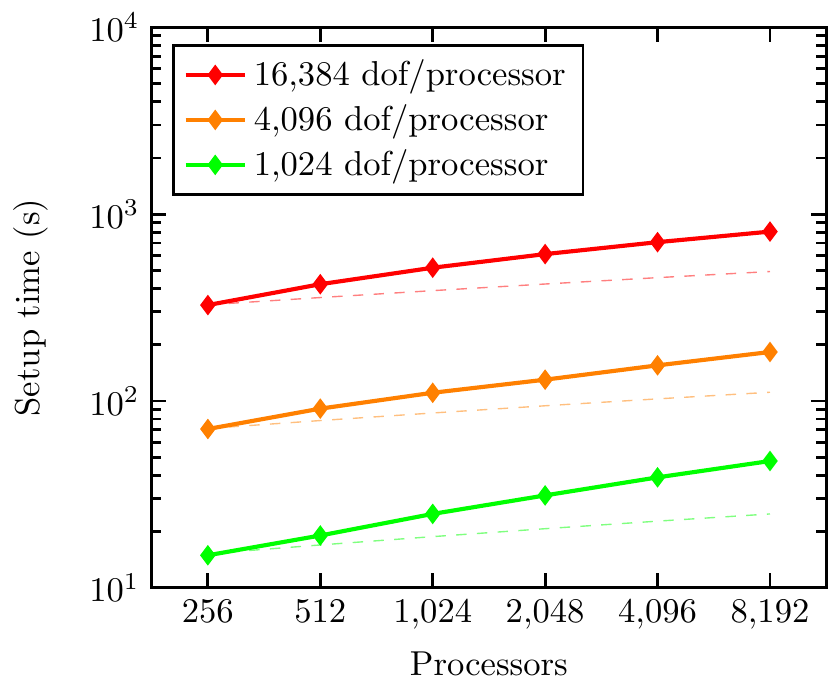}
\caption{Weak scaling of the preconditioner setup.}
\label{fig:weak1}
\end{subfigure}
\begin{subfigure}{0.4\textwidth}
\centering
\includegraphics[width=0.85\linewidth]{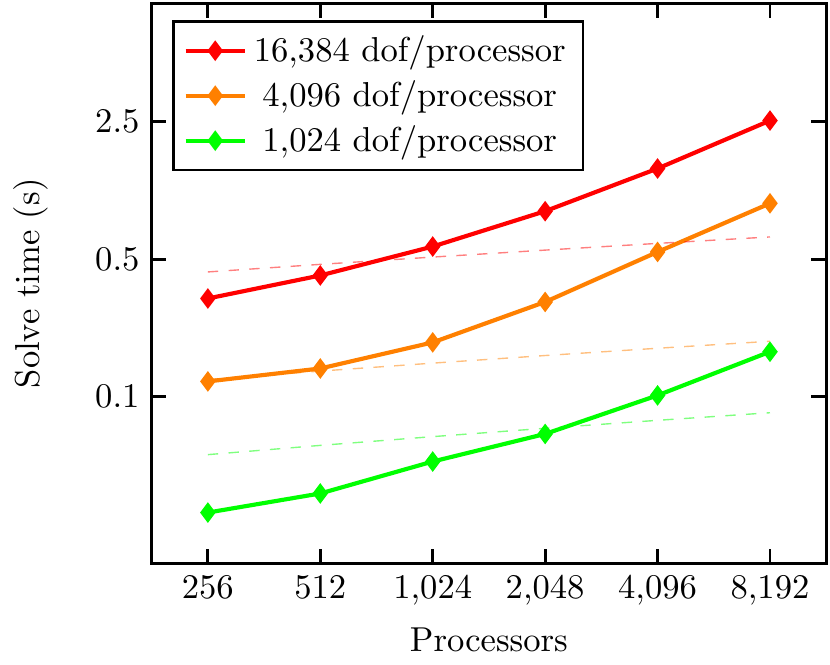}
\caption{Weak scaling of the preconditioner apply.}
\label{fig:weak2}
\end{subfigure}
\caption{Weak scalability of the ACR preconditioner application for the solution of the constant-coefficient Poisson equation.}
\end{figure}

Figures \ref{fig:weak1} and \ref{fig:weak2} depict the results of a weak scaling experiments for the ACR preconditioner fixing a different numbers of degrees of freedom per processor, along with the ideal weak scaling reference lines depicted as dashed curves considering that the estimates of setup is of $O(k^2 N \log^2 N)$ operations and the solve stage per iteration are of $O(k N \log N)$ operations.

The setup stage follows the ideal trending line as we increase the number of processors. The solve phase deviates from the ideal scaling due to the communication latency which is more noticeable due to its lower arithmetic intensity, the scalability of the Krylov method (conjugate gradient in this case), and the load imbalance in the late stages of the recursive bisection of cyclic reduction.

\section{Variable-coefficient Poisson equation}

The solution of variable-coefficient PDEs is an essential engineering problem, as the coefficient structure typically corresponds to material properties of the problem under consideration. This section documents the behavior of the ACR preconditioner from an increasingly challenging coefficient structure with up to six orders of magnitude of contrast.

The problem under consideration in this section is the symmetric positive definite discretization of the 3D variable-coefficient Poisson equation with Dirichlet boundary conditions. In particular, the second-order accurate 7-point finite-difference star stencil with harmonic average of the coefficient $\kappa(x)$ \cite{kadioglu2008comparative}:

\begin{equation}
\begin{aligned}
-\nabla \cdot \kappa(x) \nabla u = 1, \;\;\; &\mathbf{x} \in \Omega = [0,1]^3,  \;\;\; u(\mathbf{x}) = 0, \; \mathbf{x} \in \Gamma, \\
\end{aligned}
\label{poissonEquation}
\end{equation}

\subsection{Generation of random permeability fields}

The generation of random permeability field $\kappa(x)$ that closely represents a porous medium for the modeling of water or oil flow is a well-defined task on its own. The experiments in this section are based on the parallel framework for the multilevel Monte Carlo approach (MLMC) described in \cite{mohring2015uncertainty}, via the Distributed and Unified Numerics Environment DUNE \cite{bastian2008generic}. The random permeability fields are defined with covariance function of the form:

\begin{equation}
\mathcal{C}(h) = \sigma^2 \; \mbox{exp} ( - ||h||_2 / \lambda), \;\;\; h \in [0,1]^3
\end{equation}

\begin{figure}[H]
\centering
\begin{subfigure}{.225\textwidth}
\centering
\includegraphics[width=0.85\linewidth]{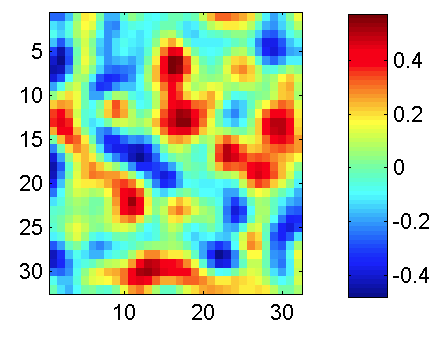}
\caption{$N = 32^3$ One order of magnitude of contrast.}
\end{subfigure}\hspace{0.1in}
\begin{subfigure}{.225\textwidth}
\centering
\includegraphics[width=0.85\linewidth]{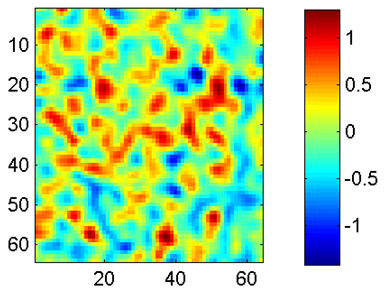}
\caption{$N = 64^3$ Two orders of magnitude of contrast.}
\end{subfigure}\hspace{0.1in}
\begin{subfigure}{.225\textwidth}
\centering
\includegraphics[width=0.85\linewidth]{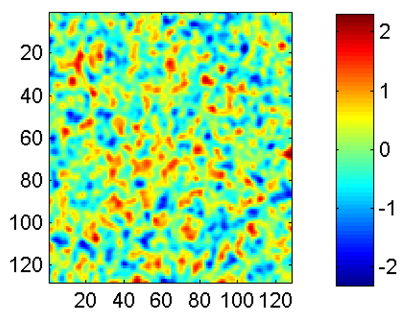}
\caption{$N = 128^3$ Four orders of magnitude of contrast.}
\end{subfigure}\hspace{0.1in}
\begin{subfigure}{.225\textwidth}
\centering
\includegraphics[width=0.85\linewidth]{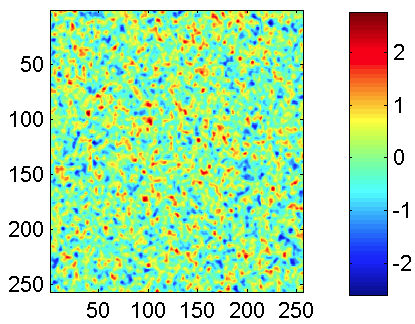}
\caption{$N = 256^3$ Six orders of magnitude of contrast.}
\end{subfigure}
\caption{Different realizations of random permeability fields  $\kappa(x)$ at different resolutions and contrast of the coefficient. Images depict the middle slice of each 3D permeability field.}
\label{figPrec:kappa}
\end{figure}

Gaussian random fields are set to a correlation length $\lambda = 3h$, where $h = \frac{1}{n-1}$ and $N=n^3$. The variance $\sigma$ is set to deliver a particular contrast in the coefficient measured in orders of magnitude. Figure \ref{figPrec:kappa} depicts four random fields realizations at different number of degrees of freedom and contrast of the coefficient.

\subsection{Tuning parameters}
\label{sec:sec1}

The main parameter that controls the accuracy of the ACR preconditioner is $\mathcal{H}_{\epsilon}$. As discussed in section \ref{sec:H}, $\mathcal{H}_{\epsilon}$ controls the accuracy of the $\mathcal{H}$-matrix approximations and their arithmetic operations. This global threshold, in turn, controls the relative accuracy of the solution for a given right-hand side.

It is expected that as we adjust $\mathcal{H}_{\epsilon}$, we can control the required number of iterations to reach convergence with a Krylov method. One could set $\mathcal{H}_{\epsilon}$ to the sought after accuracy of the solution and not require any iteration at all. However, the performance and memory requirements, although asymptotically-optimal, become impractical at high-accuracy for 3D problems.
For the ACR preconditioner, the sweet spot for achieving the fastest time to solution is not the one corresponding to the least number of iterations. It is generally the case that the inexpensive ACR preconditioners provide the fastest time to solution. There is a trade-off between the accuracy of the preconditioner and the number of Krylov iterations as numerical experiments show below.

\begin{figure}[H]
\begin{subfigure}{.5\textwidth}
\centering
\includegraphics[width=0.75\linewidth]{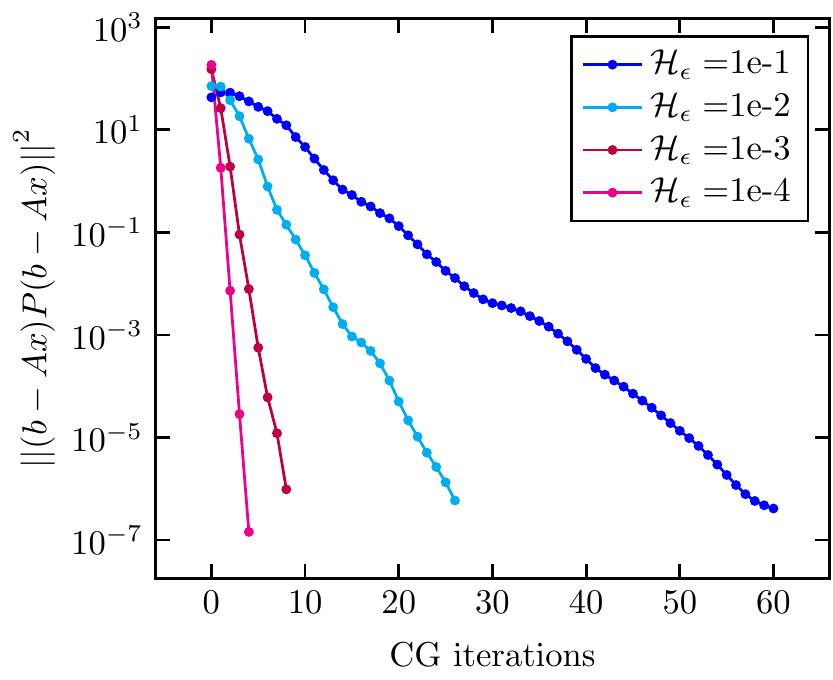}
\caption{Number of CG iterations as a function of the preconditioner accuracy $\mathcal{H}_{\epsilon}$ for the variable-coefficient Poisson equation. The preconditioner with the smallest $\mathcal{H}_{\epsilon}$ requires the least number of iterations.}
\label{figPrec:3DPoisTuning}
\end{subfigure}\hspace{0.1in}
\begin{subfigure}{.5\textwidth}
\centering
\includegraphics[width=0.75\linewidth]{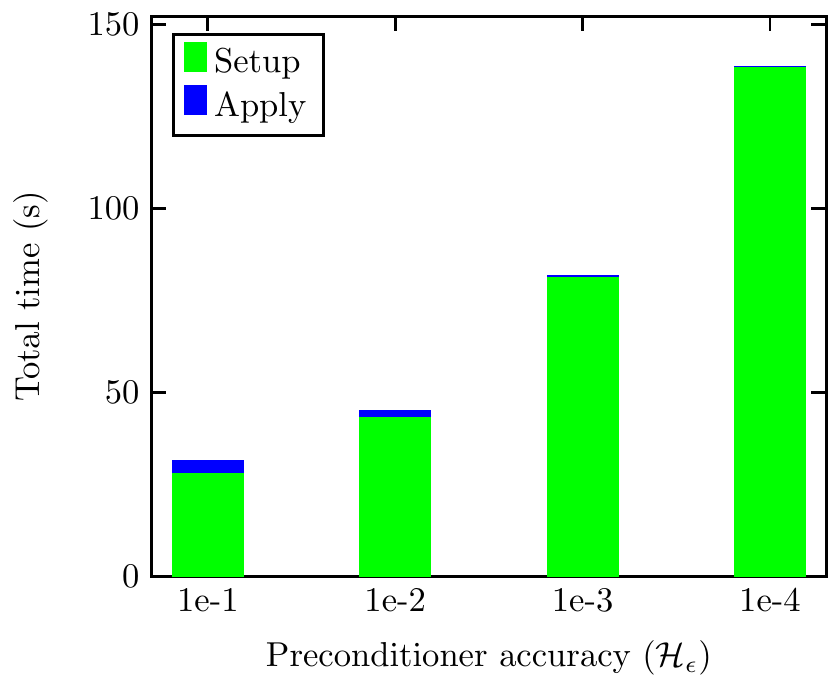}
\caption{Time requirements while refining the preconditioner accuracy $\mathcal{H}_{\epsilon}$ for the variable-coefficient Poisson equation. The preconditioner with the largest $\mathcal{H}_{\epsilon}$ delivers the best time to solution.}
\label{figPrec:3DPoisTimings}
\end{subfigure}
\caption{Number of iterations and preconditioning accuracy for the variable-coefficient Poisson equation with $N=128^3$ degrees of freedom and coefficient contrast of four orders of magnitude.}
\label{figPrec:3DPoisVar1}
\end{figure}

The effect on the required number of iterations as a function of the preconditioner accuracy $\mathcal{H}_{\epsilon}$ to solve a $N=128^3$ problem with coefficient contrast of four orders of magnitude can be seen in Figure \ref{figPrec:3DPoisTuning}. The largest $\mathcal{H}_{\epsilon}$ requires the most number of iterations, while the smallest $\mathcal{H}_{\epsilon}$ requires the least number of iterations.

\begin{figure}[H]
\begin{subfigure}{.5\textwidth}
\centering
\includegraphics[width=0.75\linewidth]{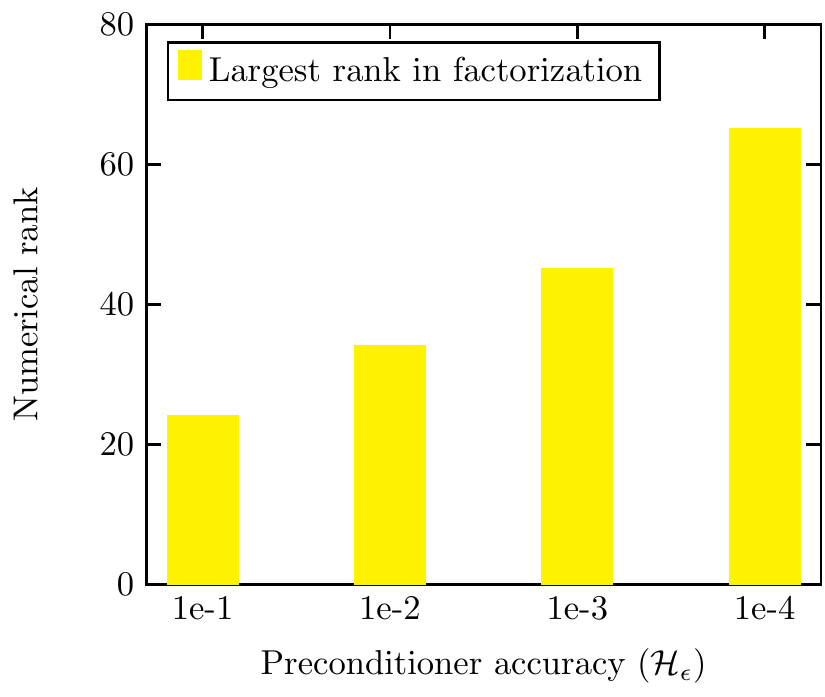}
\caption{Largest rank in the factorization at different preconditioner accuracy $\mathcal{H}_{\epsilon}$. The preconditioner with the largest $\mathcal{H}_{\epsilon}$ requires the smallest numerical rank.}
\label{figPrec:3DPoisVarRanks}
\end{subfigure}\hspace{0.1in}
\begin{subfigure}{.5\textwidth}
\centering
\includegraphics[width=0.75\linewidth]{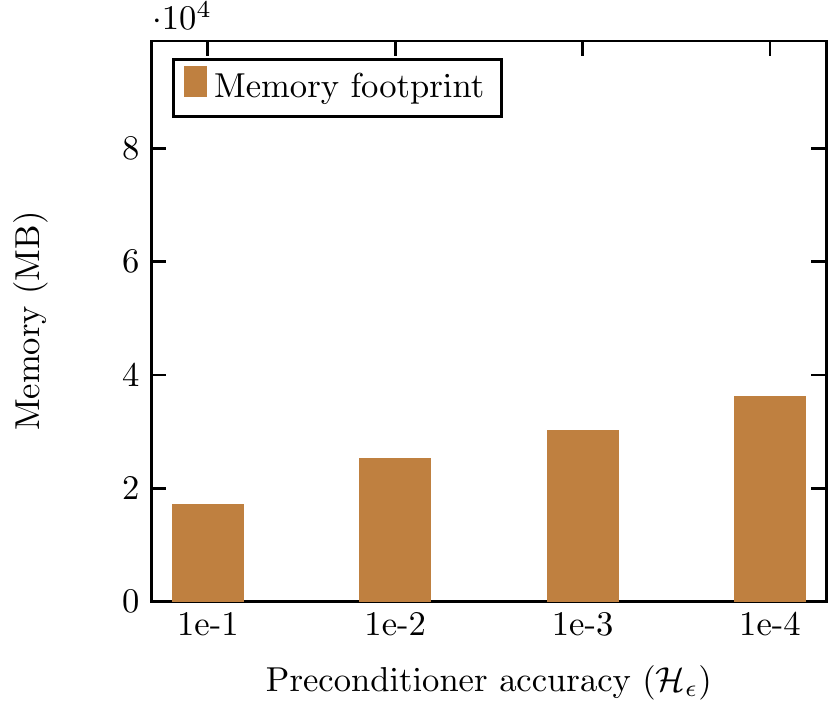}
\caption{Memory requirements at different preconditioner accuracy $\mathcal{H}_{\epsilon}$. The preconditioner with the largest $\mathcal{H}_{\epsilon}$ requires the least amount of memory.}
\label{figPrec:3DPoisVarMem}
\end{subfigure}
\caption{Effect of the preconditioner accuracy $\mathcal{H}_{\epsilon}$ for the variable-coefficient Poisson equation with $N=128^3$ degrees of freedom and coefficient contrast of four orders of magnitude.}
\label{figPrec:3DPoisVar2}
\end{figure}

As can be seen from Figure \ref{figPrec:3DPoisTimings}, even though setting a large $\mathcal{H}_{\epsilon}$ required the greatest number of CG iterations, this is the recommended value of $\mathcal{H}_{\epsilon}$ to optimize for time to solution in our current implementation. Although there are more iterations than with a smaller $\mathcal{H}_{\epsilon}$, the application of the preconditioner is fastest at large $\mathcal{H}_{\epsilon}$ since the ranks are the smallest, see Figure \ref{figPrec:3DPoisVarRanks}. Figure \ref{figPrec:3DPoisVarMem} depicts how $\mathcal{H}_{\epsilon}$ directly determines the memory footprint of the preconditioner, and shows why it is desirable to set $\mathcal{H}_{\epsilon}$ as large as possible to also optimize for memory requirements. Larger values of the preconditioner accuracy could deliver better time to solution, although at the expense of more synchronizing iterations, which might be undesirable for some applications.

\subsection{Sensitivity with respect to high contrast coefficient}

As the problem difficulty increases, i.e. the contrast of the coefficient sharpens, there are cases for which the most economical preconditioner (e.g. $\mathcal{H}_{\epsilon}$=1e-1) might not reach convergence within an acceptable number of iterations, see Table \ref{tableNotConvergence}.

\begin{table}[H]
\centering
\begin{tabular}{|c|c|c|}
\hline
$N$ & $\mathcal{H}_{\epsilon}$ & CG Iterations \\ \hline
$32^3$                   & 1e-1 & 27         \\ \hline
$64^3$                   & 1e-1 & 51         \\ \hline
$128^3$                  & 1e-1 & 95         \\ \hline
\multirow{2}{*}{$256^3$} & 1e-1 & 100+        \\ \cline{2-3} 
& 1e-2 & 73         \\ \hline
\end{tabular}
\caption{Number of iterations required by CG for the variable-coefficient Poisson equation with coefficient contrast of six orders of magnitude. The most economical preconditioner for the hardest problem did not reach convergence within 100 iterations, thus requiring a more accurate version of the preconditioner to reach convergence.}
\label{tableNotConvergence}
\end{table}

\begin{figure}[H]
\centering
\includegraphics[width=0.4\linewidth]{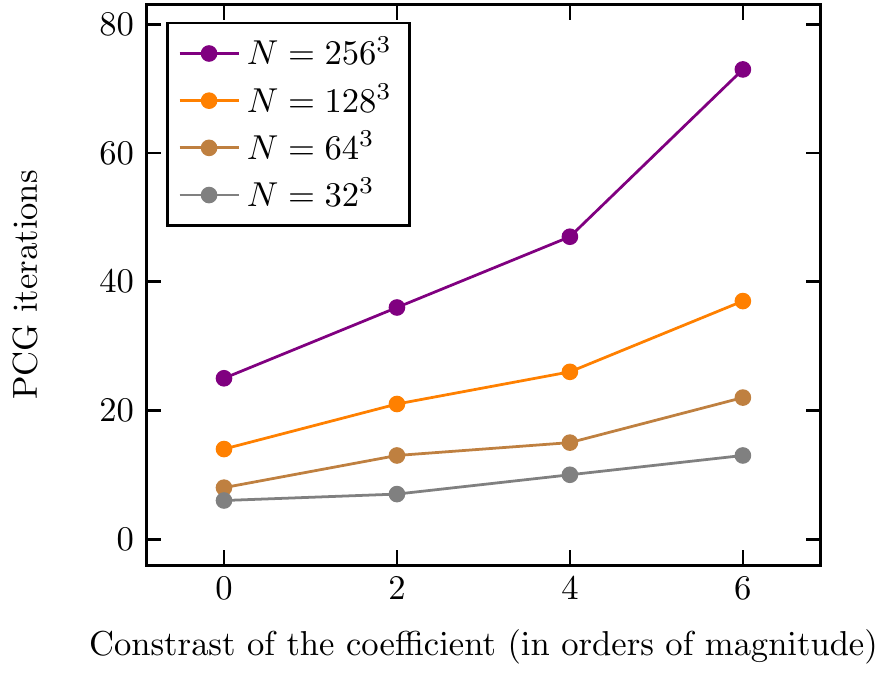}
\caption{Required number of iterations for an ACR preconditioner accuracy of $\mathcal{H}_{\epsilon}$=1e-2 as both the problem size and the contrast of the coefficient increases. A larger number of iterations is necessary as the contrast of the coefficient increases.}
\label{figPrec:3DPoisIncreaseContrast}
\end{figure}

Therefore, to reach convergence, a more accurate preconditioner is necessary. Figure \ref{figPrec:3DPoisIncreaseContrast} shows the required number of iterations to achieve convergence for a preconditioner with accuracy $\mathcal{H}_{\epsilon}$=1e-2 at increasing problem size and contrast of the coefficient. For comparison, the baseline case (zero contrast) depicts a constant-coefficient Poisson equation with $\kappa(x)$=1.

One of the established methods for the solution of these type of problems is algebraic multigrid. As the contrast of the coefficient changes, Table \ref{table:comparisonsPoisson} shows the number of iterations required by the conjugate gradient method (CG) without preconditioning, and the preconditioned conjugate gradient method using algebraic multigrid (AMG) as a preconditioner, accessed here via the hypre library \cite{Briggs00, hypre_web_page}. To compare at a similar number of iterations, the ACR preconditioner was tuned to require less than ten iterations, in this case at $\mathcal{H}_{\epsilon}$=1e-4. Experiments show that at moderate contrast of the coefficients and using 512 cores, the solve time of the ACR preconditioner is comparable to the AMG preconditioner, albeit at a large enough number of right-hand sides so that the setup time of ACR gets amortized since the preconditioner can be reused. Another dimension of comparison is a traditional direct solver, such as the LU factorization, here accessed through the SuperLU DIST package \cite{li2003superlu_dist}. For instance, for the problem of four orders of magnitude of contrast, the LU factors are computed in 27.16 seconds and require 3.5E4 MB of memory, whereas the ACR preconditioner at $\mathcal{H}_{\epsilon}$=1e-1 is computed in 27.71 seconds, but it only requires 1.7E4 MB of memory to store its factors.

\begin{table}[H]
\centering
\begin{tabular}{|c|c|c|c|c|c|c|}
\hline
\multirow{2}{*}{\begin{tabular}[c]{@{}c@{}}Coefficient\\ contrast\end{tabular}} & \multicolumn{2}{c|}{CG + No Prec.} & \multicolumn{2}{c|}{CG + AMG} & \multicolumn{2}{c|}{CG + ACR} \\ \cline{2-7} 
        & Iterations         & Solve         & Iterations       & Solve      & Iterations       & Solve      \\ \hline
0                                                                               & 257                & 0.10          & 6                & 0.53       & 3                & 0.30       \\ \hline
2                                                                               & 975                & 0.24          & 6                & 0.35       & 4                & 0.34       \\ \hline
4                                                                               & 2210               & 0.44          & 6                & 0.28       & 4                & 0.37       \\ \hline
6                                                                               & 6968               & 1.35          & 7                & 0.26       & 7                & 0.59       \\ \hline
\end{tabular}
\caption{Number of iterations and solve time for the solution of a sequence of Poisson problems $N=128^3$ with a variable coefficient at increasing order of magnitude of contrast of the coefficient. Methods under consideration include algebraic multigrid (AMG), and accelerated cyclic reduction (ACR).}
\label{table:comparisonsPoisson}
\end{table}

\subsection{Operation count and memory footprint}

\begin{figure}[H]
\begin{subfigure}{.5\textwidth}
\centering
\includegraphics[width=0.75\linewidth]{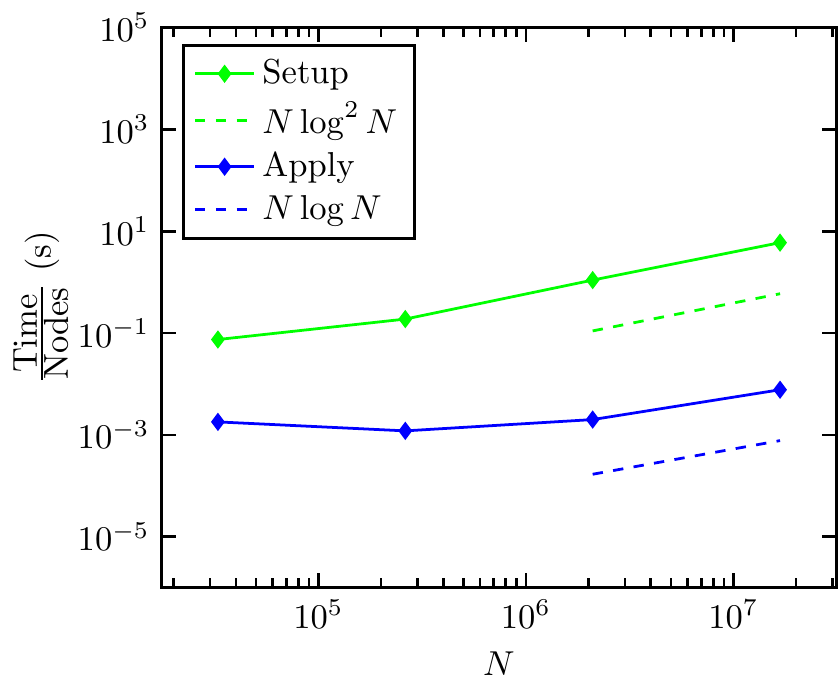}
\caption{Comparison of the preconditioner setup and application with their corresponding theoretical estimates.}
\label{figPrec:3DPoisSetupApply}
\end{subfigure}\hspace{0.1in}
\begin{subfigure}{.5\textwidth}
\centering
\includegraphics[width=0.75\linewidth]{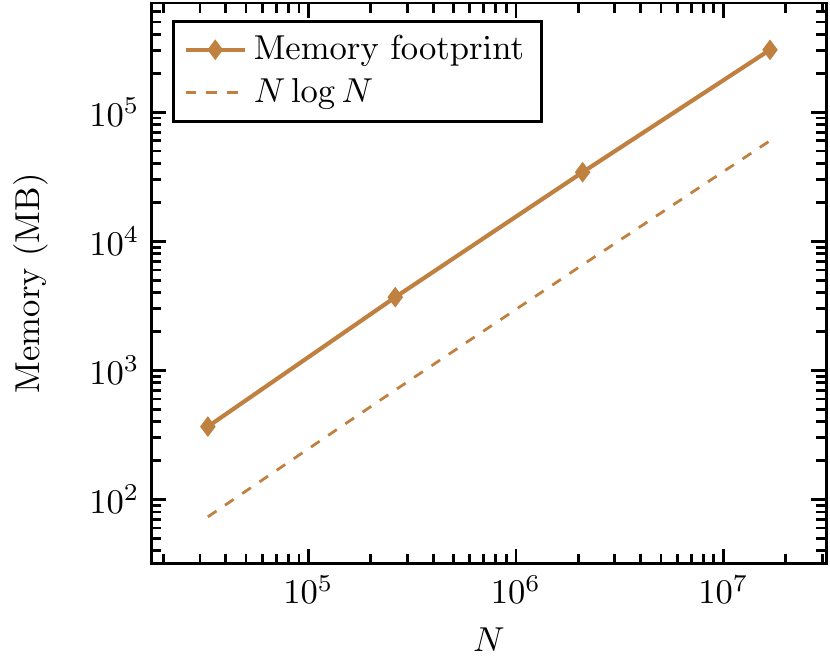}
\caption{Comparison of the preconditioner memory footprint with its theoretical estimate.}
\label{figPrec:3DPoisMem}
\end{subfigure}
\caption{Measured performance and memory footprint for the solution of an increasingly larger variable-coefficient Poisson equation with a random field of four orders of magnitude of contrast in the coefficient. The preconditioner accuracy for this experiments is set to $\mathcal{H}_{\epsilon}$ = 1e-1.}
\label{figPrec:scaleN}
\end{figure}

The complexity estimates of the number of operations in the setup and application phases of the preconditioner are bounded by $O(k^2N \log^2 N)$ and $O(kN \log N)$ respectively, while its memory footprint is bounded by $O(kN \log N)$; where $N$ is the number of degrees of freedom and $k$ is the numerical rank of the approximation. To demonstrate that these estimates hold for a variable-coefficient problem, Figure \ref{figPrec:scaleN} shows the behavior of the preconditioner in terms of operations count and memory footprint as we increase the number of degrees of freedom $N$ for the variable-coefficient Poisson equation with a coefficient of four orders of magnitude of contrast.

The vertical axis of Figure \ref{figPrec:3DPoisSetupApply}, normalized by the number of processors used in each case, reports the measured performance of the setup and application phases of the preconditioner while comparing it with their theoretical complexity. Figure \ref{figPrec:3DPoisMem} reports the total memory requirements as the problem size increases and also compares it with the theoretical complexity demonstrating a fair agreement.

\section{Convection-diffusion equation with recirculating flow}

In this section we show the effectiveness of the ACR preconditioner on the convection-diffusion equation with a variable and recirculating flow $b(\mathbf{x})$, i.e. a flow with vanishing normal velocities at the boundary.

\begin{equation}
\begin{aligned}
&-\nabla \cdot \kappa(x) \nabla u  + \alpha b(\mathbf{x}) \cdot \nabla u = f(\mathbf{x}), \;\;\; \mathbf{x} \in \Omega = [0,1]^3, \\
&b(\mathbf{x}) = \begin{bmatrix}
\sin(a\;2\pi x) \sin(a\;2\pi (1/8 + y)) +  \sin(a\;2\pi (1/8 + z)) \sin(a\;2\pi x) \\
\cos(a\;2\pi x) \cos(a\;2\pi (1/8 + y)) +  \cos(a\;2\pi (1/8 + y)) \cos(a\;2\pi z) \\
\cos(a\;2\pi x) \cos(a\;2\pi (1/8 + z)) +  \sin(a\;2\pi (1/8 + y)) \sin(a\;2\pi z)
\end{bmatrix},\\
&b_x+b_y+b_z=0.
\end{aligned}
\label{eqPrec:codi}
\end{equation}

The equation is discretized with a 7-point upwind finite difference scheme, which leads to a nonsymmetric linear system. When the convection term dominates, $\alpha >1$, this equation is known to be challenging for classical iterative solvers.

\subsection{Tuning parameters}
\label{label:TuningParameters}

\begin{figure}[H]
\begin{subfigure}{.5\textwidth}
\centering
\includegraphics[width=0.75\linewidth]{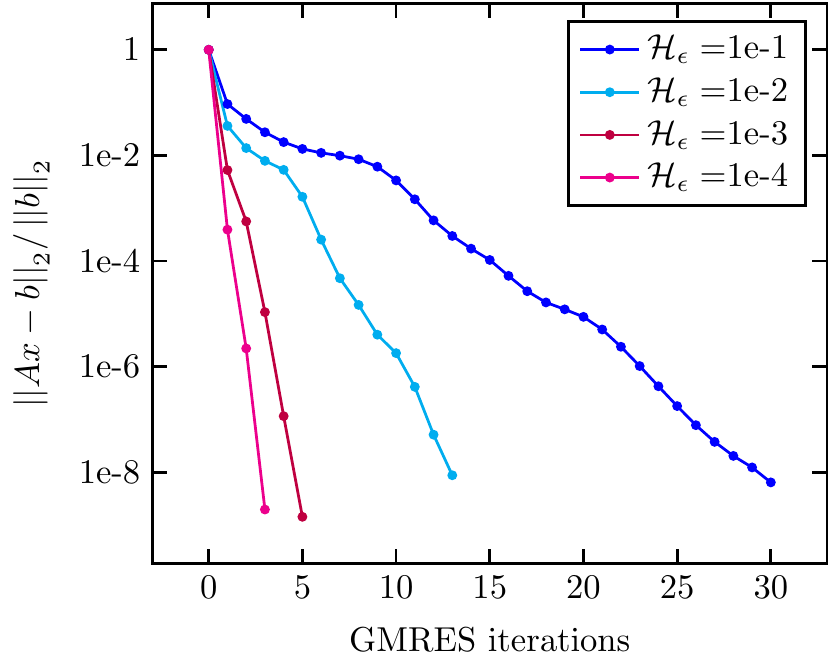}
\caption{Number of iterations as a function of the preconditioner accuracy $\mathcal{H}_{\epsilon}$. As $\mathcal{H}_{\epsilon}$ decreases, the preconditioner requires fewer iterations.}
\label{figPrec:3DCodiConvergence}
\end{subfigure}\hspace{0.1in}
\begin{subfigure}{.5\textwidth}
\centering
\includegraphics[width=0.75\linewidth]{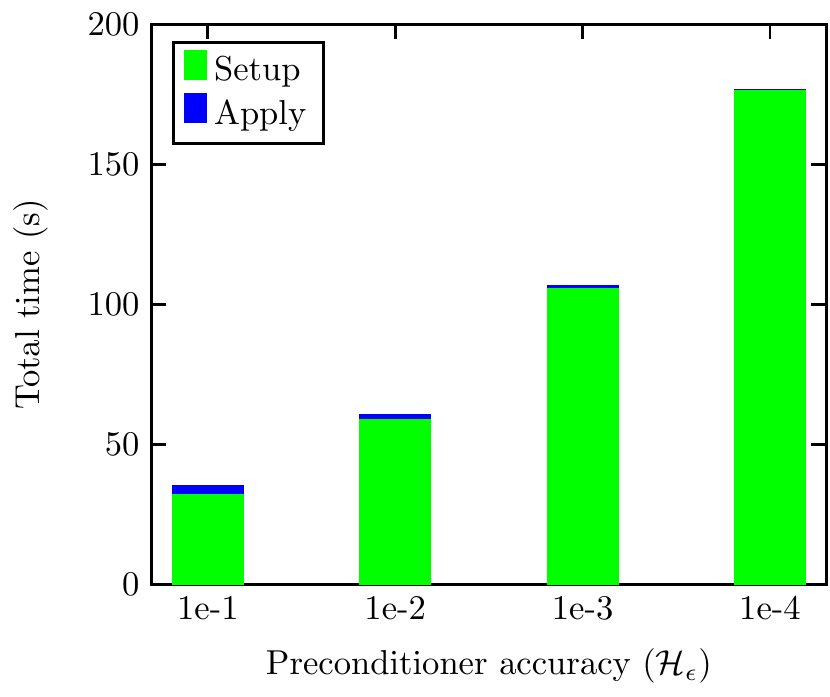}
\caption{Time requirements while refining the preconditioner accuracy $\mathcal{H}_{\epsilon}$. The largest $\mathcal{H}_{\epsilon}$ delivers the best time to solution.}
\label{figPrec:3DCodiTimings}
\end{subfigure}
\caption{This experiment depicts a convection-diffusion problem with recirculating flow with eight vortices, $\alpha=8$, discretized with $N=128^3$ degrees of freedom.}
\end{figure}

In a regime of convection dominance, Figure \ref{figPrec:3DCodiConvergence} shows how the ACR preconditioner can control the number of GMRES iterations by tuning the preconditioner accuracy $\mathcal{H}_{\epsilon}$.

\begin{figure}[H]
\begin{subfigure}{.5\textwidth}
\centering
\includegraphics[width=0.75\linewidth]{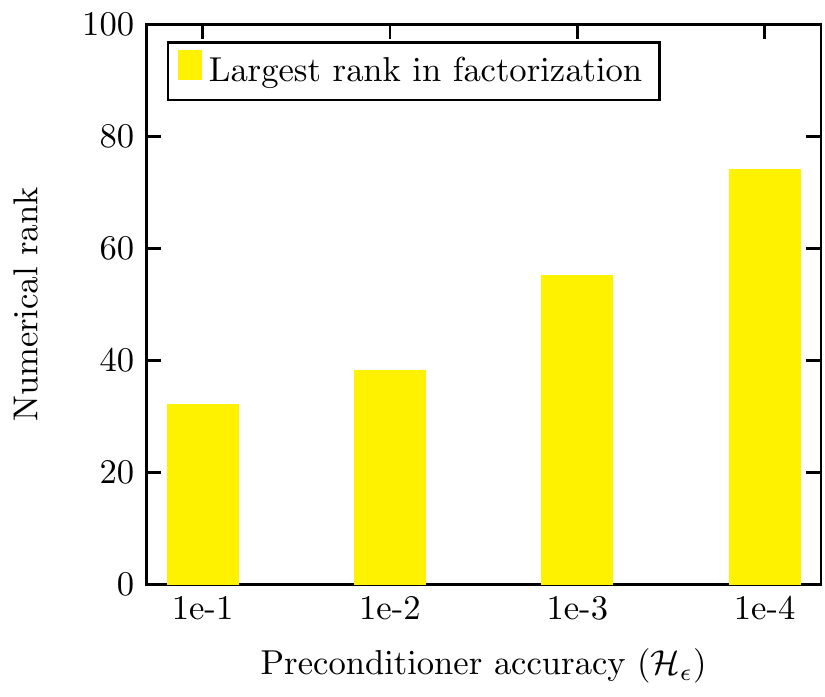}
\caption{Largest rank in factorization at different $\mathcal{H}_{\epsilon}$ for the convection-diffusion equation. The preconditioner with the largest $\mathcal{H}_{\epsilon}$ features the lowest numerical ranks.}
\label{figPrec:3DCoDiRanks}
\end{subfigure}\hspace{0.1in}
\begin{subfigure}{.5\textwidth}
\centering
\includegraphics[width=0.75\linewidth]{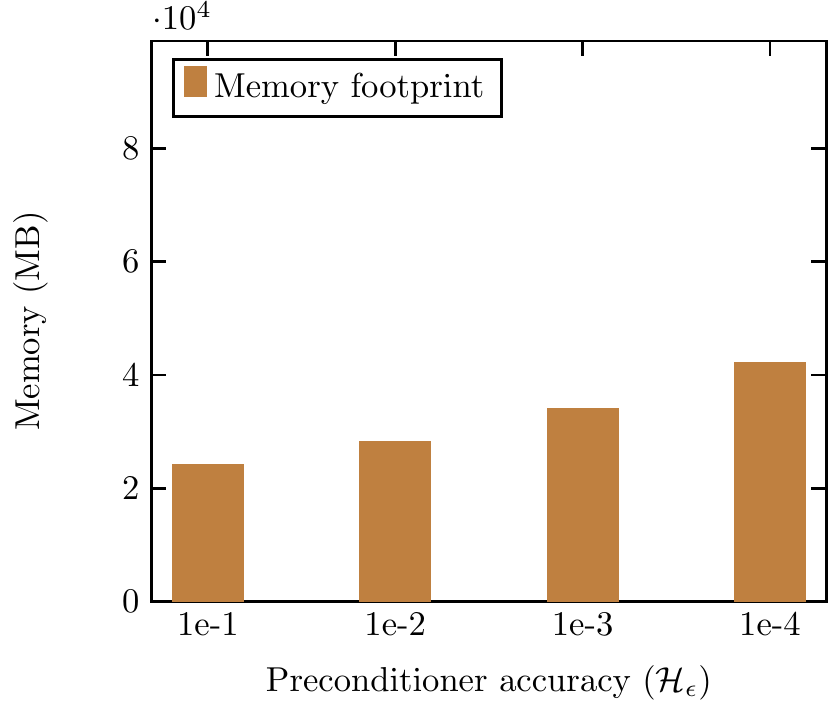}
\caption{Memory requirements while refining the preconditioner accuracy $\mathcal{H}_{\epsilon}$. The largest $\mathcal{H}_{\epsilon}$ delivers the preconditioner with the least memory footprint.}
\label{figPrec:3DCoDiMem}
\end{subfigure}
\caption{Effect on the preconditioner accuracy $\mathcal{H}_{\epsilon}$ for a convection-diffusion problem with recirculating flow with eight vortices, $\alpha=8$, and discretized with $N=128^3$ degrees of freedom.}
\end{figure}

Regarding absolute time to solution, experiments with our latest implementation show that the preconditioner with the largest $\mathcal{H}_{\epsilon}$ led to the best time to solution, albeit with the most iterations, as Figure \ref{figPrec:3DCodiTimings} shows. As a result, this preconditioner configuration featured the lowest numerical rank, as shown in Figure \ref{figPrec:3DCoDiRanks}, enabling a fast application at each iteration. Furthermore, the fastest preconditioner had the least memory requirements, as shown in Figure \ref{figPrec:3DCoDiMem}.

\subsection{Sensitivity with respect to vortex wavenumber}

\begin{figure}[H]
\centering
\begin{subfigure}{0.24\textwidth}
\centering
\includegraphics[width=0.85\linewidth]{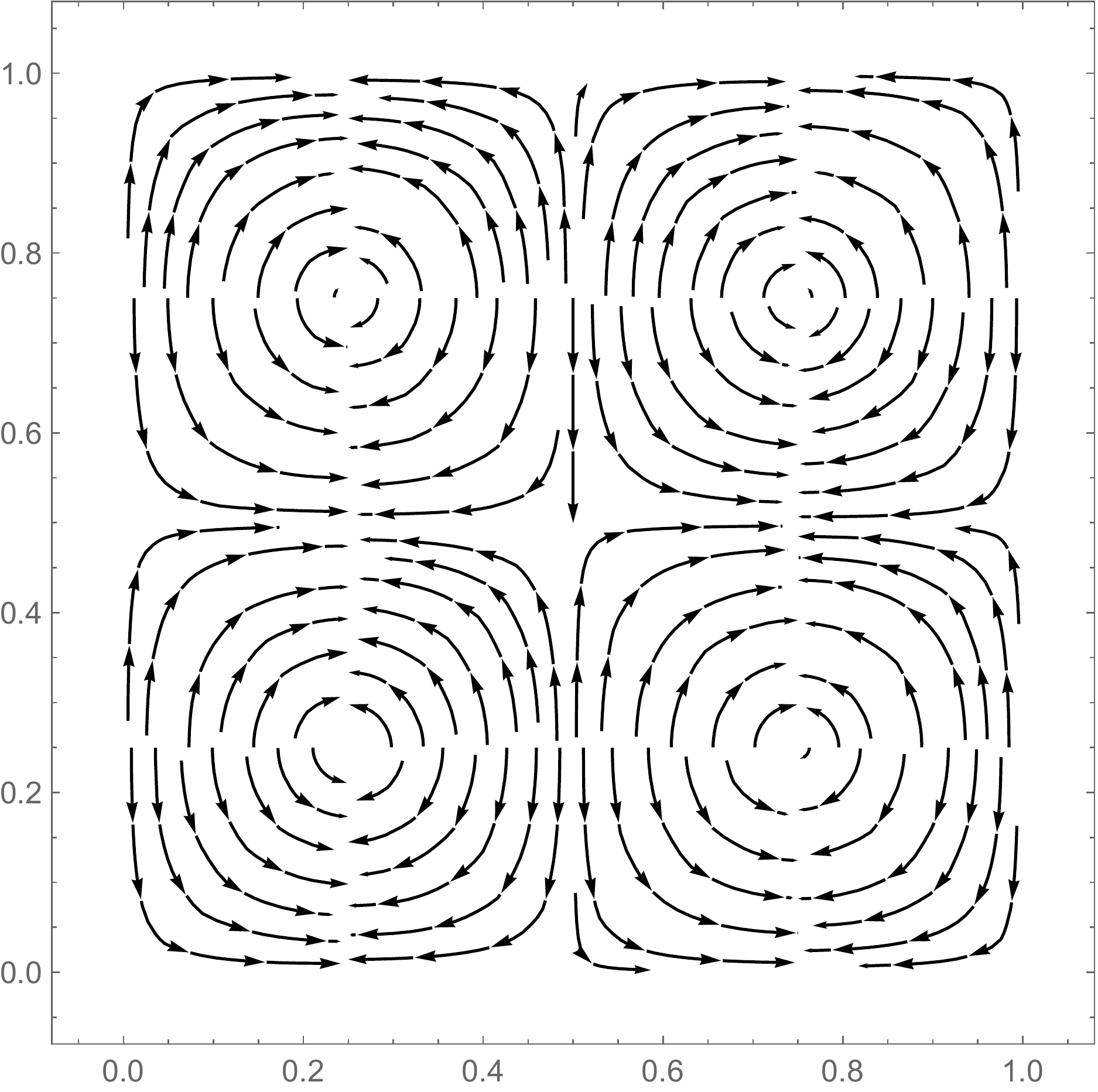}
\caption{Two vortices.}
\end{subfigure}
\begin{subfigure}{0.24\textwidth}
\centering
\includegraphics[width=0.85\linewidth]{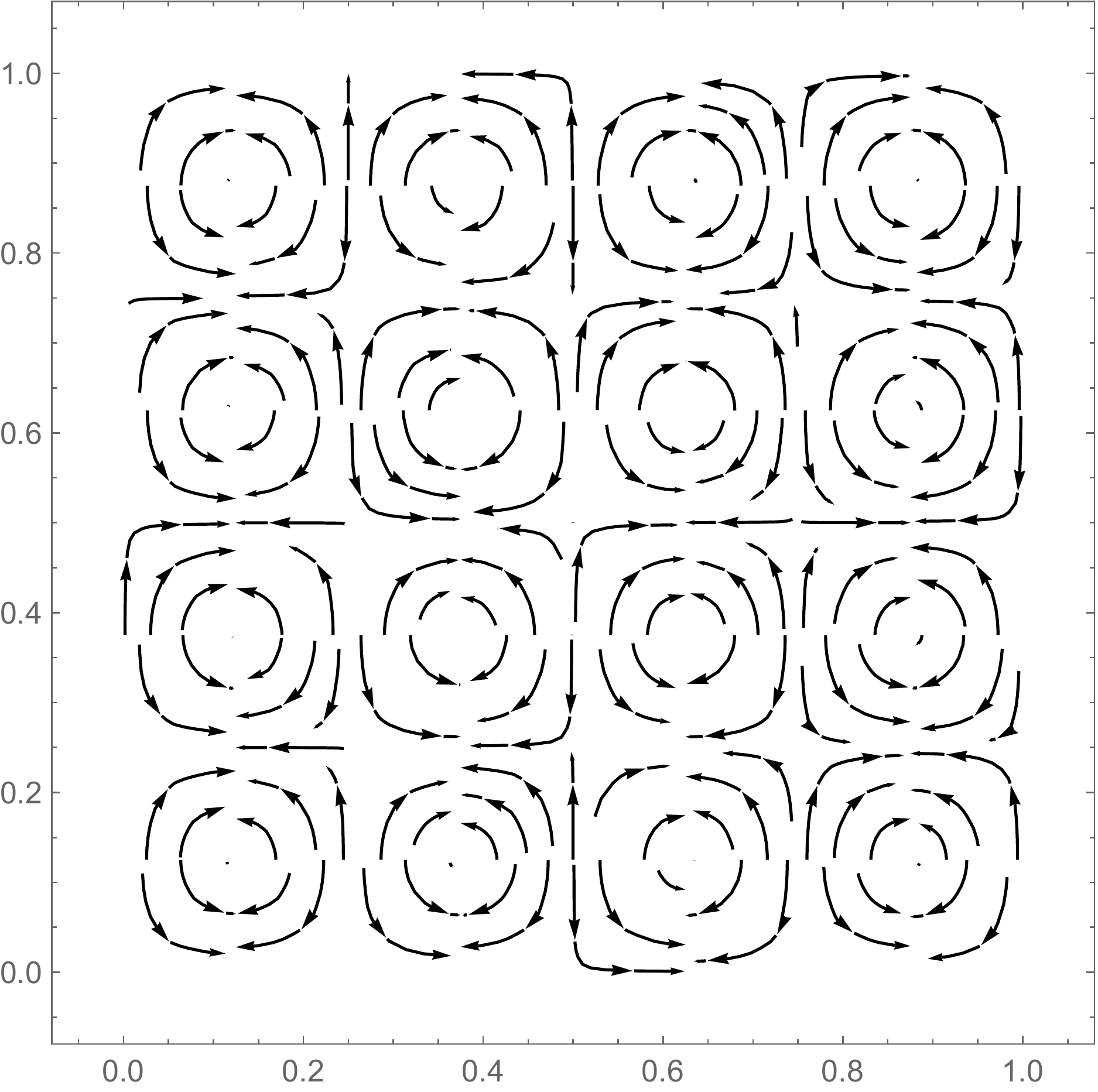}
\caption{Four vortices.}
\end{subfigure}
\begin{subfigure}{0.24\textwidth}
\centering
\includegraphics[width=0.85\linewidth]{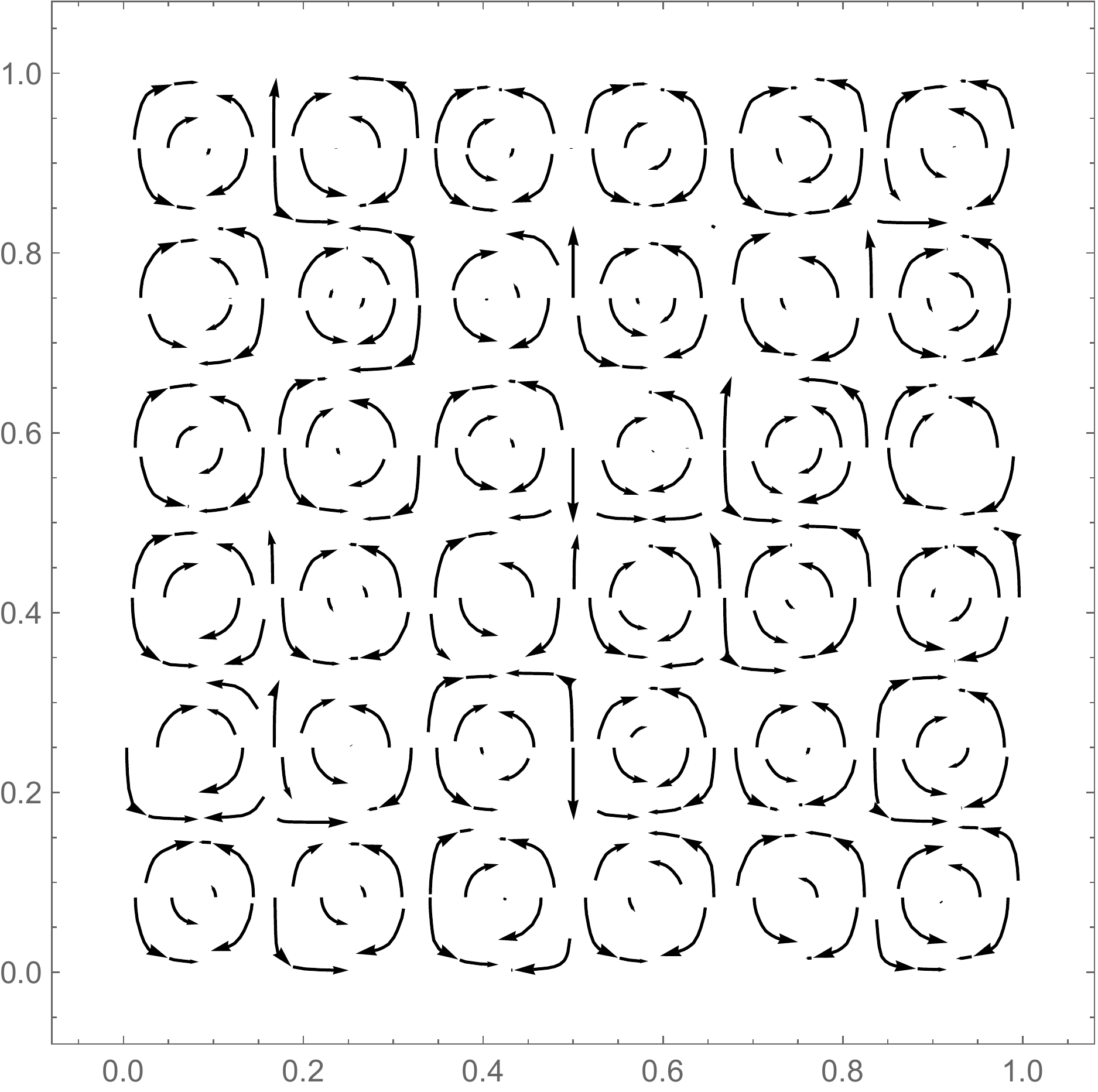}
\caption{Six vortices.}
\end{subfigure}
\begin{subfigure}{0.24\textwidth}
\centering
\includegraphics[width=0.85\linewidth]{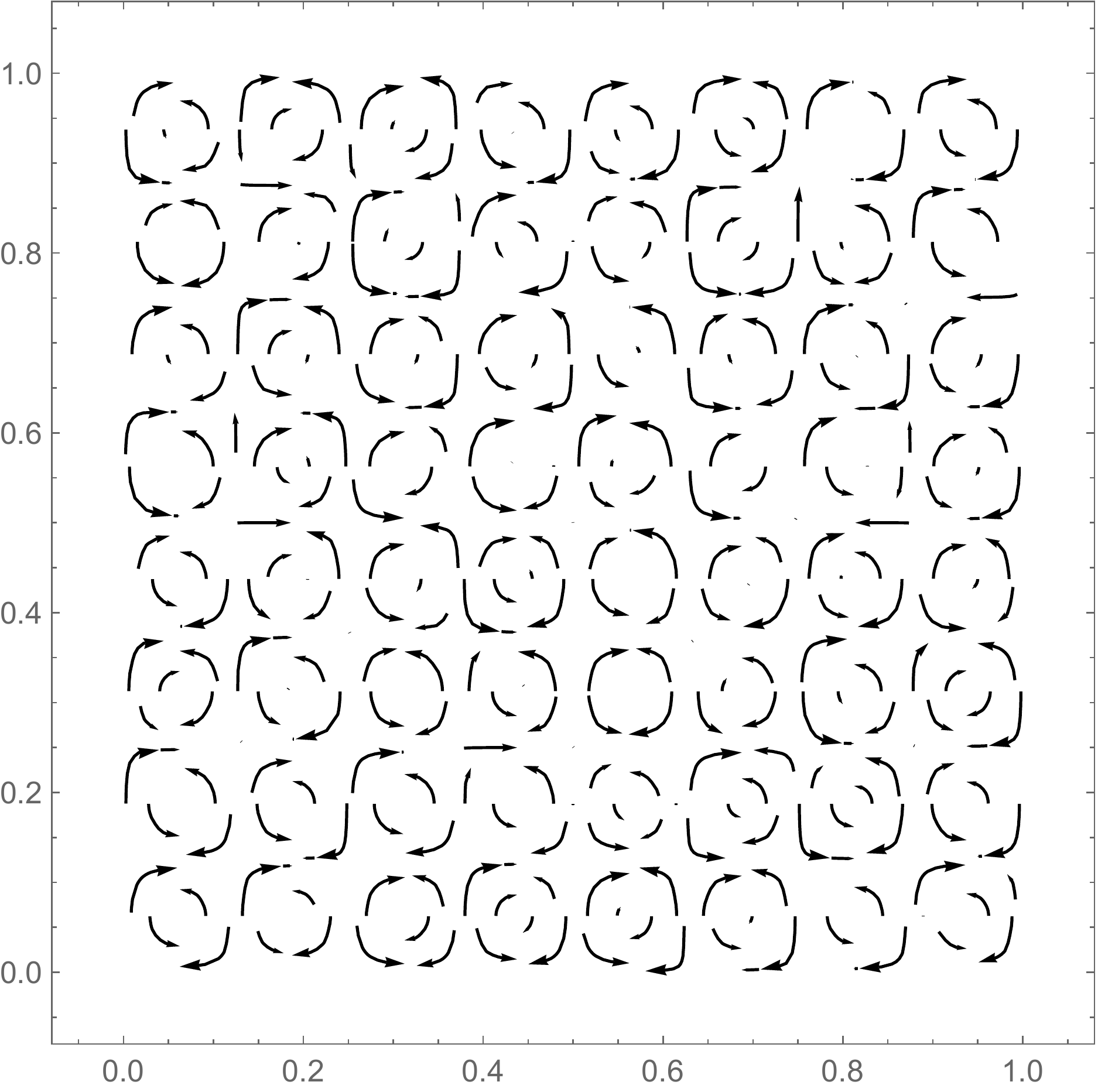}
\caption{Eight vortices.}
\end{subfigure}
\caption{Increasing number of vortices in the flow $b(\mathbf{x})$.}
\label{figPrec:MultipleVortices}
\end{figure}

\begin{figure}[H]
\centering
\includegraphics[width=0.4\linewidth]{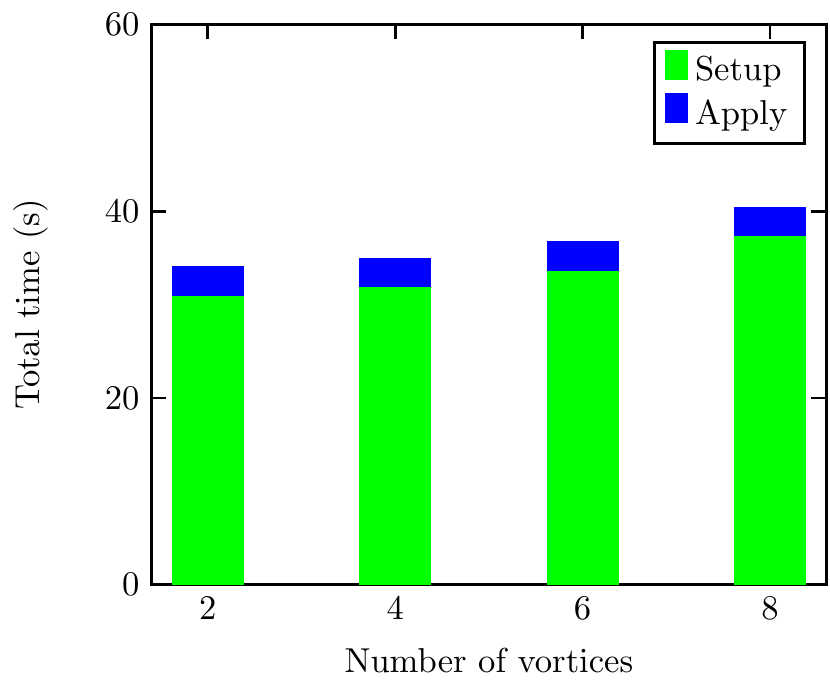}
\caption{Time distribution of the preconditioner as the number of vortices in $b(\mathbf{x})$ increases. Increasing the number of vortices had a minor effect on the effectiveness of the preconditioner.}
\label{figPrec:CoDiMultipleVortices}
\end{figure}

Consider an increasing number of vortices in the flow $b(\mathbf{x})$, as Figure \ref{figPrec:MultipleVortices} shows. At the corners, and in center of each vortex, there are saddle points which are known to be challenging for multigrid methods to resolve \cite{gupta2000high}. Figure \ref{figPrec:CoDiMultipleVortices} demonstrates that the ACR preconditioner remains robust as the number of vortices increases.

\subsection{Sensitivity with respect to convection dominance}

\begin{figure}[H]
\begin{subfigure}{0.5\textwidth}
\centering
\includegraphics[width=0.75\linewidth]{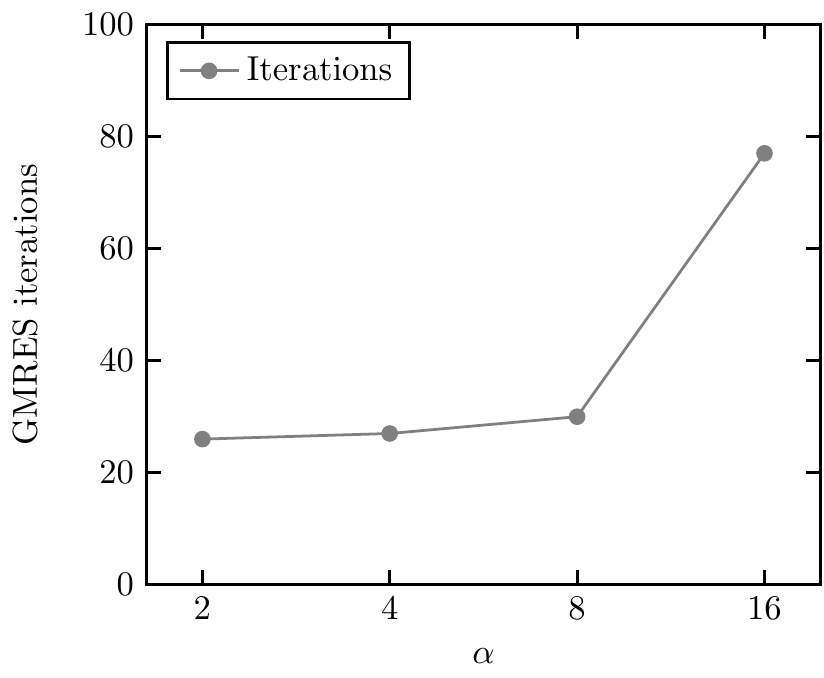}
\caption{Number of iterations as the convection term gains dominance. An increase in the dominance of the convection term requires mores iterations.}
\label{figPrec:3DCoDiDominantConvectionIterations}
\end{subfigure}\hspace{0.1in}
\begin{subfigure}{.5\textwidth}
\centering
\includegraphics[width=0.75\linewidth]{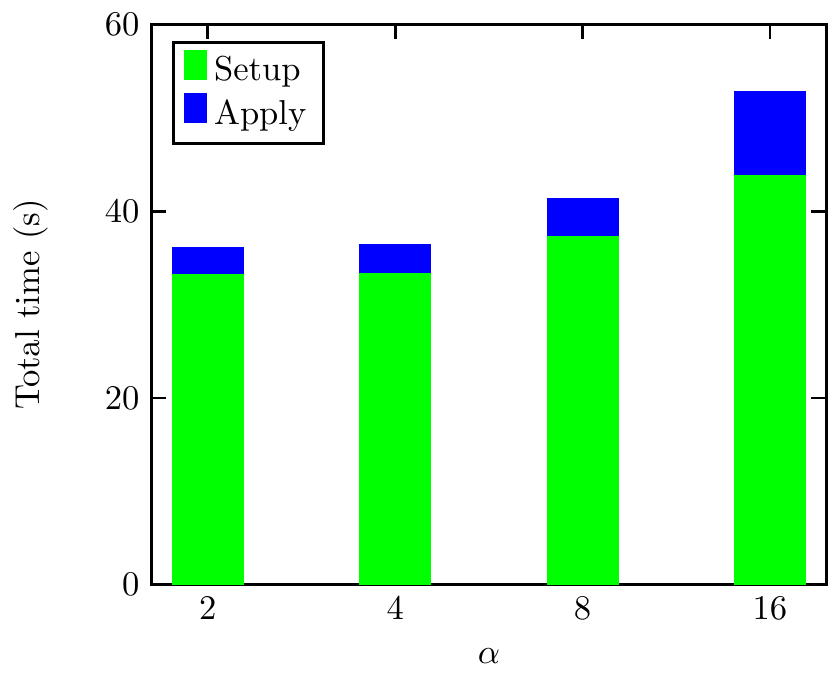}
\caption{Time requirements as the convection term gains dominance. The overall time to solution has a moderate increase.}
\label{figPrec:3DCoDiDominantConvectionTime}
\end{subfigure}
\caption{Effect on the preconditioner accuracy $\mathcal{H}_{\epsilon}$ for the convection-diffusion equation with recirculating flow discretized with $N=128^3$ degrees of freedom as the convective becomes more significant than the diffusion term.}
\label{figPrec:3DCoDiDominantConvection}
\end{figure}

Consider a fixed-accuracy ACR preconditioner with $\mathcal{H}_{\epsilon}=$1e-1, and an increasingly convection dominated problem, achieved by gradually increasing $\alpha$ in Equation \ref{eqPrec:codi}. As expected, Figure \ref{figPrec:3DCoDiDominantConvectionIterations} shows that a low-accuracy preconditioner requires more iterations as the convective term dominates. Furthermore, given that the accuracy of the preconditioner is fixed, there is a noticeable effect on the application phase of the preconditioner, which is proportional to the number of iterations. A graphical representation of such behavior can be seen in Figure \ref{figPrec:3DCoDiDominantConvectionTime}. Evidently, as shown in the section \ref{label:TuningParameters}, it is possible to control, and decrease, the number of iterations by building a more accurate preconditioner. But the key point here is that the ACR preconditioner in combination with GMRES is demonstrated to be robust for convection dominated problems.

As a matter of comparison, Table \ref{table:comparisonsCoDi} shows the number of iterations of GMRES without preconditioner, and GMRES in combination with the incomplete LU factorization. In this case, the problem with the strongest convection dominance did not reach a solution with the incomplete LU factorization as a preconditioner of GMRES. In terms of memory, the ACR preconditioner for the problem with the strongest convection dominance required 30.76 seconds and 2.3E4 MB of memory with $\mathcal{H}_{\epsilon}$=1e-1, whereas SuperLU DIST required 44.59 seconds and 3.7E4 MB of memory.

\begin{table}[H]
\centering
\begin{tabular}{|c|c|c|c|c|c|c|}
\hline
\multirow{2}{*}{$\alpha$} & \multicolumn{2}{c|}{GMRES + No Prec.} & \multicolumn{2}{c|}{GMRES + ILU(0)} & \multicolumn{2}{c|}{GMRES + ACR} \\ \cline{2-7} 
& Iterations              & Solve                & Iterations     & Solve         & Iterations        & Solve        \\ \hline
0                         & 1,132                & 0.38           & 201                 & 0.20          & 26                & 2.93         \\ \hline
2                         & 1,242                & 0.51           & 226                 & 0.22          & 27                & 3.10         \\ \hline
4                         & 1,765                & 0.60           & 368                 & 0.30          & 30                & 4.02         \\ \hline
6                         & 100,000+             & -              & 100,000+            & -             & 77                & 8.94         \\ \hline
\end{tabular}
\caption{Number of iterations and solve time for the solution of a sequence of convection-diffusion problems with $N=128^3$ and increasingly convection dominance. Methods under consideration include the incomplete LU factorization (ILU), and accelerated cyclic reduction (ACR).}
\label{table:comparisonsCoDi}
\end{table}

\subsection{Operation count and memory footprint}

\begin{figure}[H]
\begin{subfigure}{.5\textwidth}
\centering
\includegraphics[width=0.75\linewidth]{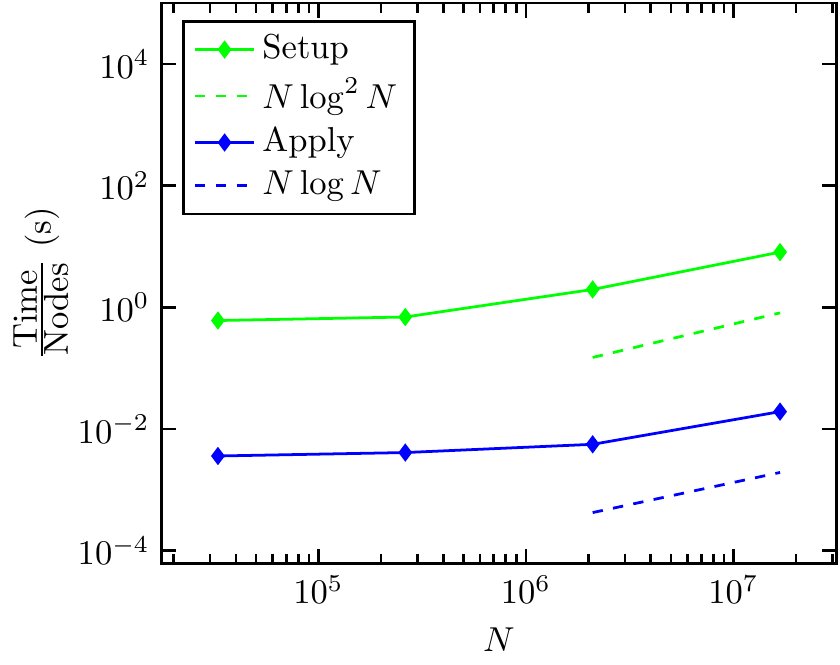}
\caption{Comparison of the preconditioner setup and application with their corresponding theoretical estimates.}
\label{figPrec:3DCoDiComplexityTime}
\end{subfigure}\hspace{0.1in}
\begin{subfigure}{.5\textwidth}
\centering
\includegraphics[width=0.75\linewidth]{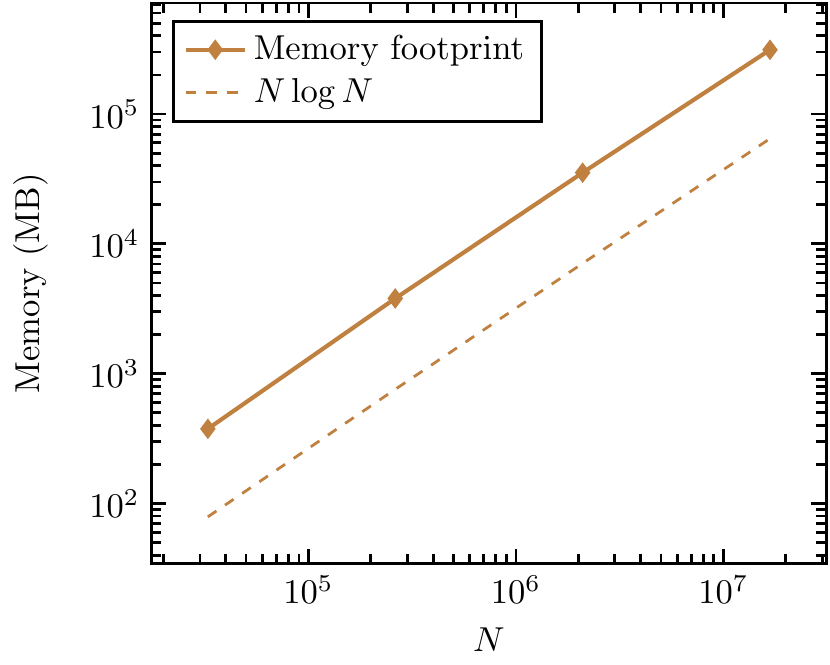}
\caption{Comparison of the preconditioner memory footprint with its theoretical estimate.}
\label{figPrec:3DCoDiComplexityMemory}
\end{subfigure}
\caption{Measured performance and memory footprint for the solution of the convection-diffusion equation with recirculating flow.}
\label{figPrec:3DCoDiComplexity}
\end{figure}

Figure \ref{figPrec:3DCoDiComplexity} presents a comparison between the measured performance and memory requirements of the preconditioner, and their corresponding theoretical complexity estimates for the convection-diffusion problem described in Equation \ref{eqPrec:codi}, as the problem size increases.

The vertical axis of Figure \ref{figPrec:3DCoDiComplexityTime}, normalized with the number of compute nodes used in each case, reports the measured performance of the setup and application phases of the preconditioner while demonstrating a fair agreement with the asymptotic complexity estimates for the large-scale experiments.

Figure \ref{figPrec:3DCoDiComplexityMemory} reports the total memory requirements as the problem size increases and also compares it with its corresponding theoretical complexity demonstrating a fair agreement across all experiments.

\section{Indefinite Helmholtz equation in heterogeneous media}
The numerical solution of the indefinite Helmholtz equation offers one of the greatest challenges for iterative and direct solvers at large-scale \cite{ernst2012difficult}. There is a significant interest in the development of optimal methods as several engineering applications use the Helmholtz equation to model time-harmonic propagation of acoustic waves. Inversion techniques based on full-waveform inversion (FWI) for instance, involve heterogeneous velocity models and the solution of multiple right-hand sides at a wide range of frequencies. Therefore, the introduction of an efficient forward solver directly contributes to expanding the limits of what can be modeled computationally.

Consider the indefinite Helmholtz equation in a variable velocity field $c(\mathbf{x})$, at frequency $f$, and Dirichlet boundary conditions in the unit cube:

\begin{equation}
\begin{aligned}
- \nabla^2u - \frac{{(2 \pi f)}^2}{c(\mathbf{x})^2}u &= f( \mathbf{x} ), \; \Omega = [0,1]^3, \; \mathbf{x} \in \Gamma, \\
c(\mathbf{x}) & = 1.25(1-0.4 e^{-32( |x-0.5|^2+|y-0.5|^2)} ) \\
u( \mathbf{x} ) & = \sin(\pi x)\sin(\pi y)\sin(\pi z)
\end{aligned}
\label{eq:HelmVarVel}
\end{equation}

\begin{figure}[H]
\centering
\includegraphics[width=0.4\linewidth]{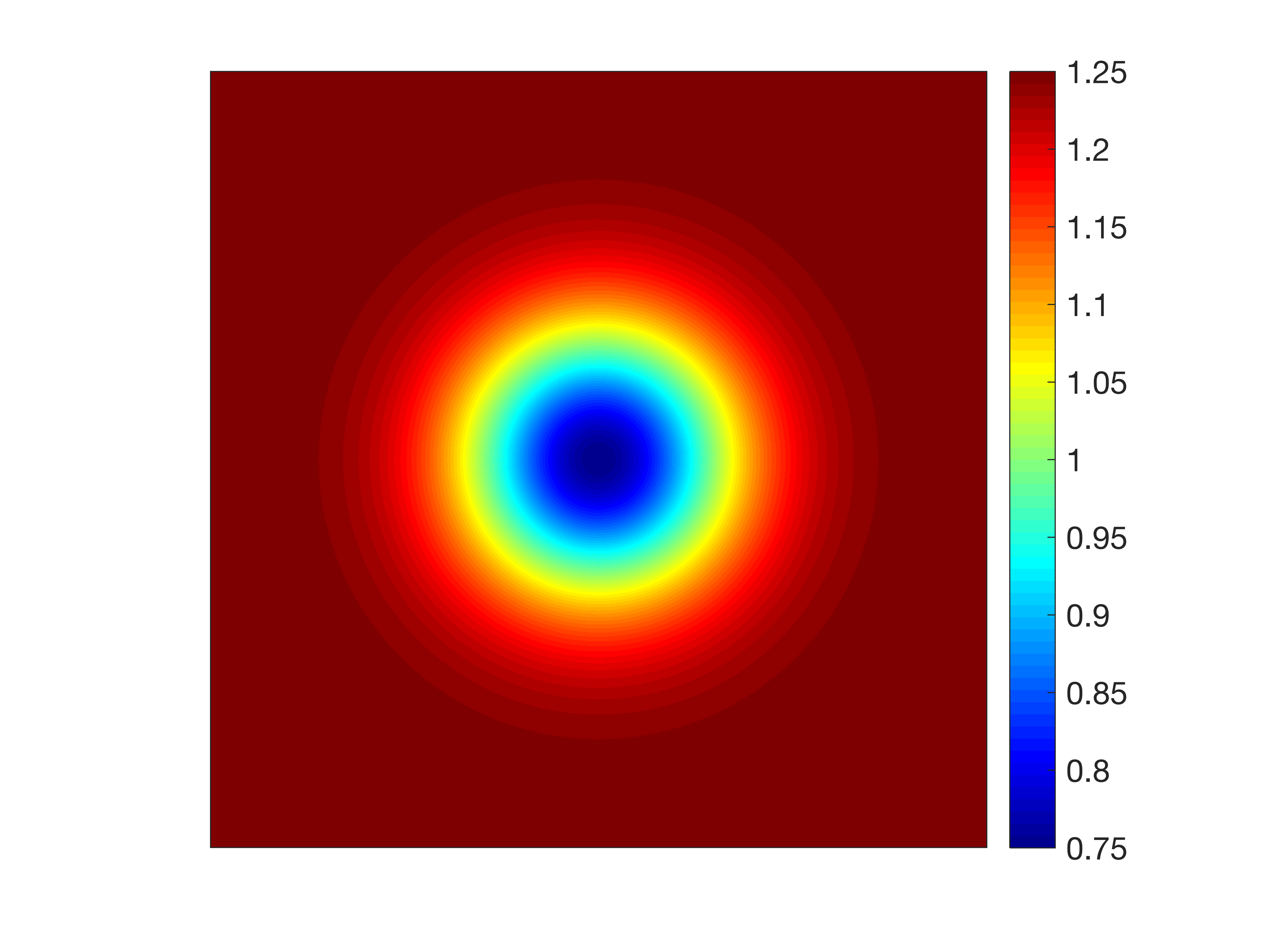}
\caption{Wave velocity field $c(\mathbf{x})$. The image depicts the middle slice of the 3D wave velocity field.}
\label{figPrec:3DHelmVelocityField}
\end{figure}

The velocity field models a waveguide over the unit cube as proposed in \cite{poulsonSweeping2013}, and depicted in \ref{figPrec:3DHelmVelocityField}. The forcing term $f( \mathbf{x} )$ is adjusted to satisfy the proposed exact solution $u( \mathbf{x} )$. The equation is discretized with the 27-point trilinear finite element scheme on hexahedra with the software library PetIGA \cite{PetIGA2013}. Since the linear system arising from the discretization is indefinite in the high-frequency Helmholtz case, we use ACR to accelerate the convergence of GMRES.

\subsection{Tuning parameters}
\label{sec:HelmRanks}

\begin{figure}[H]
\centering
\includegraphics[width=0.4\linewidth]{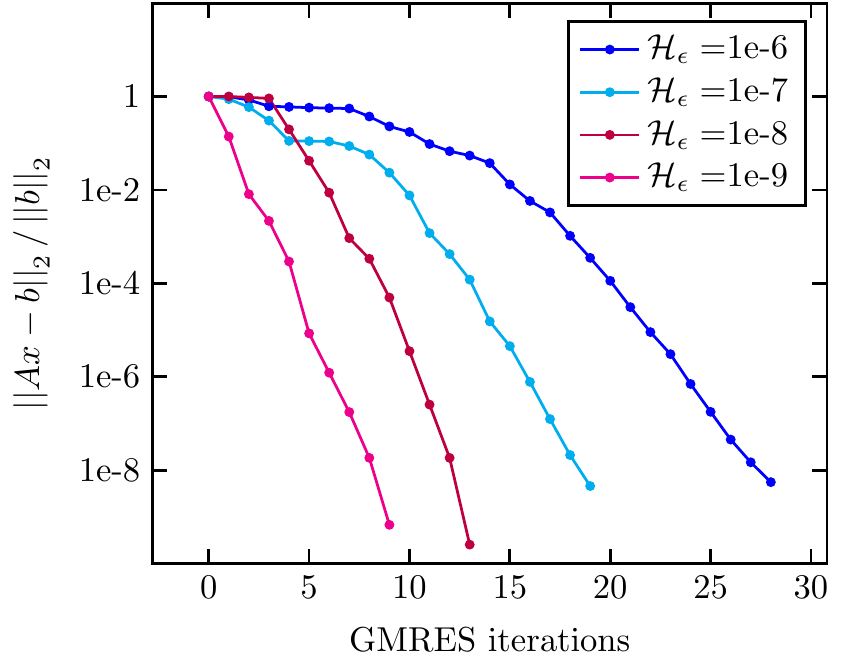}
\caption{Number of iterations as a function of the preconditioner accuracy $\mathcal{H}_{\epsilon}$. As $\mathcal{H}_{\epsilon}$ decreases, the preconditioner requires fewer iterations.}
\label{figPrec:3DHelmVarTuning}
\end{figure}

We illustrate the effectiveness of the ACR preconditioner on a moderately high-frequency Helmholtz problem as described in Equation \ref{eq:HelmVarVel}, discretized with $N=128^3$ degrees of freedom and 12 points per wavelength. As Figure \ref{figPrec:3DHelmVarTuning} shows, we can control the number of iterations that GMRES requires to reach convergence by adjusting the accuracy of the preconditioner $\mathcal{H}_{\epsilon}$. Notice that the preconditioner accuracy $\mathcal{H}_{\epsilon}$ is smaller than what was chosen for diffusive problems. The need of higher relative accuracy is due to the fact that the Helmholtz equation, in the high-frequency regime, has off-diagonal block ranks that asymptotically grow with problem size ($k \sim \mathcal{O}(n)$). This theoretical estimate is reported in the literature \cite{chandrasekaran2010}. Evidently, rank growth impacts hierarchical-matrix based solvers. Nonetheless, the complexity estimates of the ACR preconditioner are still lower than traditional exact sparse factorizations, as demonstrated in section \ref{sec:HelmEstimates}.

\begin{figure}[H]
\begin{subfigure}[H]{.5\textwidth}
\centering
\includegraphics[width=0.75\linewidth]{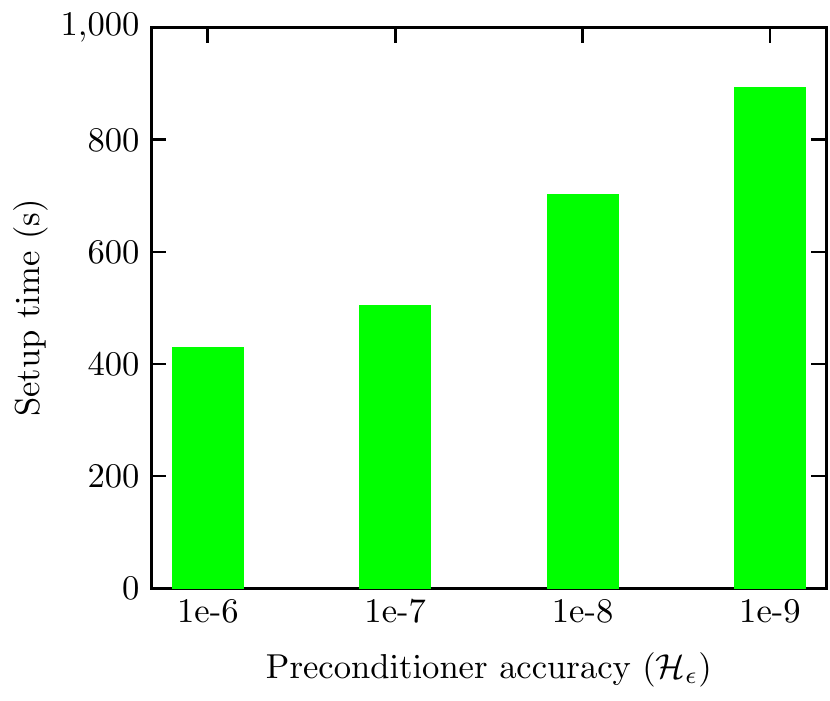}
\caption{Setup phase at increasing accuracy of the preconditioner.}
\label{figPrec:3DHelmVarTimings_setup}
\end{subfigure}\hspace{0.1in}
\begin{subfigure}[H]{.5\textwidth}
\centering
\includegraphics[width=0.75\linewidth]{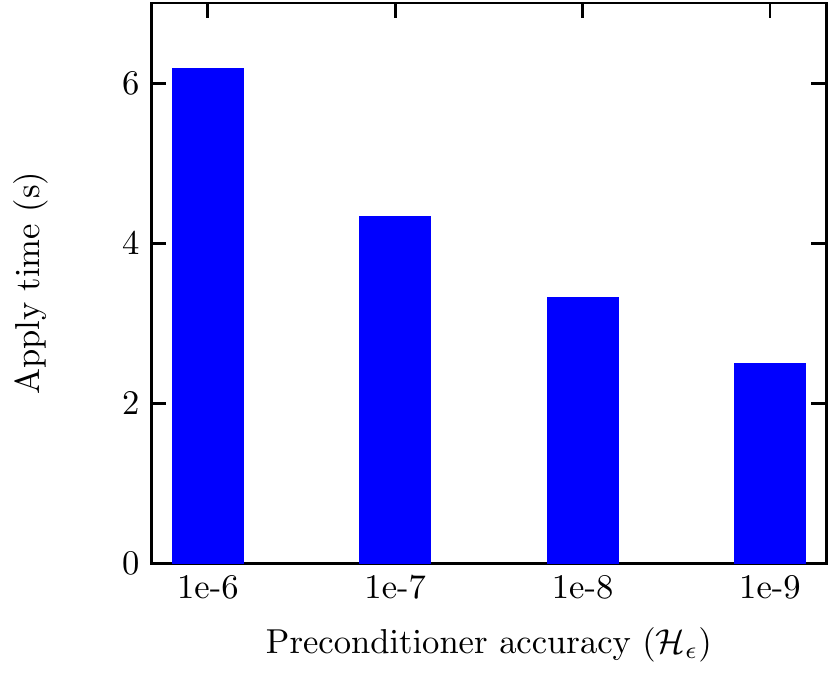}
\caption{Application phase at increasing accuracy of the preconditioner.}
\label{figPrec:3DHelmVarTimings_apply}
\end{subfigure}\hspace{0.1in}
\caption{Time requirements while refining the preconditioner accuracy $\mathcal{H}_{\epsilon}$. The loosest $\mathcal{H}_{\epsilon}$ delivers the fastest time to solution for a single right-hand side, whereas the tightest $\mathcal{H}_{\epsilon}$ delivers the best time to solution for a large number of right-hand sides, since the preconditioner setup is computed only once.}
\label{figPrec:3DHelmVarTimings}
\end{figure}

Even though the timings of the preconditioner shows an increase in the setup time as compared to diffusive problems, it still features an economical solve stage (Figure \ref{figPrec:3DHelmVarTimings_apply}). As mentioned in the introduction of this section, for inverse problems which require the solution of a large number of right-hand sides (typically up to a few thousands), the setup phase (Figure \ref{figPrec:3DHelmVarTimings_setup}) is typically regarded as an off-line phase that gets amortized if the solve stage is relatively fast, which is the case for the ACR preconditioner.

\begin{figure}[H]
\begin{subfigure}{0.5\textwidth}
\centering
\includegraphics[width=0.75\linewidth]{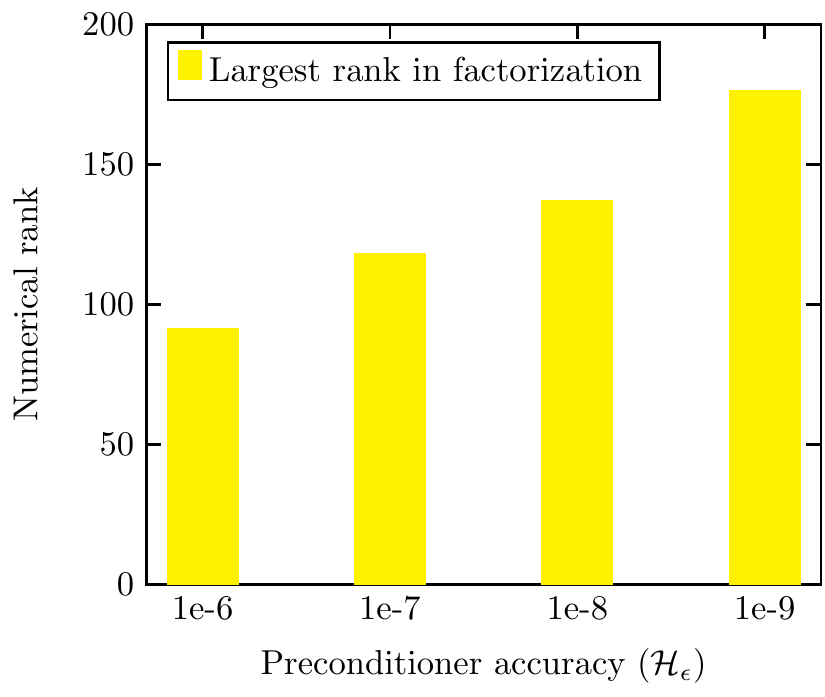}
\caption{Largest rank in factorization. Factorizations with smaller ranks lead to more iterations, but less time to solution and memory footprint.}
\label{figPrec:3DHelmVarRanks}
\end{subfigure}\hspace{0.1in}
\begin{subfigure}{0.5\textwidth}
\centering
\includegraphics[width=0.75\linewidth]{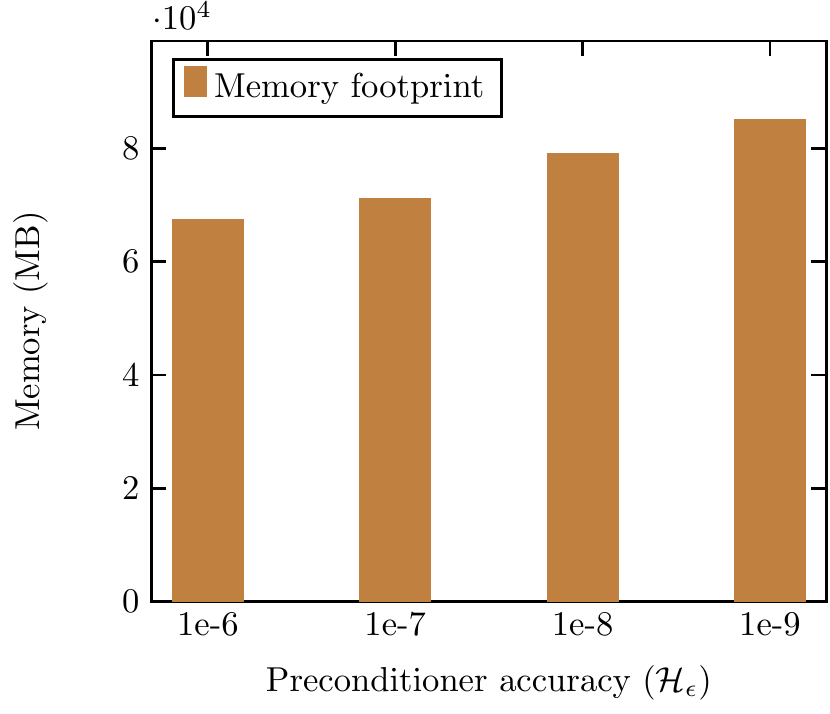}
\caption{Memory requirements while refining the preconditioner accuracy $\mathcal{H}_{\epsilon}$. The most economical preconditioner regarding memory footprint is delivered with the largest $\mathcal{H}_{\epsilon}$.}
\label{figPrec:3DHelmVarMem}
\end{subfigure}
\caption{Effect on the preconditioner accuracy $\mathcal{H}_{\epsilon}$ for the high-frequency Helmholtz equation in a heterogeneous medium discretized with $N=128^3$ degrees of freedom and 12 points per wavelength.}
\label{figPrec:3DHelmVarRanksAndMem}
\end{figure}

The growth in the setup phase as the accuracy of the preconditioner is tightened is due to increased numerical ranks, as shown in Figure \ref{figPrec:3DHelmVarRanks}. Rank growth has a direct impact on the memory footprint of the preconditioner, as shown in Figure \ref{figPrec:3DHelmVarMem}. Once more, the preconditioner with the loosest $\mathcal{H}_{\epsilon}$, i.e. the lowest numerical rank, is the preconditioner of choice to optimize for both memory and performance.

\subsection{Low to high frequency Helmholtz regimes}

Consider a sequence of Helmholtz problems, as described in Equation \ref{eq:HelmVarVel}, at increasing frequency. If the frequency is set to $f=0$ Hz, the zeroth-order term vanishes, and we are left with a constant-coefficient Poisson problem. At the other end of the spectrum, a frequency of $f=8$ Hz corresponds to a moderately high-frequency Helmholtz problem at 12 points per wavelength, as the problem featured in the previous section. Table \ref{table:AdjustHelmMultiFreq} shows the preconditioner accuracy chosen to require a maximum of 20 GMRES iterations to reach convergence.

\begin{table}[H]
\centering
\begin{tabular}{|c|c|c|}
\hline
$f$ & \begin{tabular}[c]{@{}c@{}}Points per\\ wavelength\end{tabular} & $\mathcal{H}_{\epsilon}$ \\ \hline
0 & -   & 1e-1 \\ \hline
2 & 48  & 1e-3 \\ \hline
4 & 24  & 1e-4 \\ \hline
8 & 12  & 1e-6 \\ \hline
\end{tabular}
\caption{Tuning of the preconditioner to require at most 20 GMRES iterations for a sequence of Helmholtz problems at increasing frequencies. The problem with $f=0$ represents a constant-coefficient Poisson problem, while $f=8$ represents a moderately high-frequency Helmholtz problem.}
\label{table:AdjustHelmMultiFreq}
\end{table}

\begin{figure}[H]
\begin{subfigure}{.5\textwidth}
\centering
\includegraphics[width=0.75\linewidth]{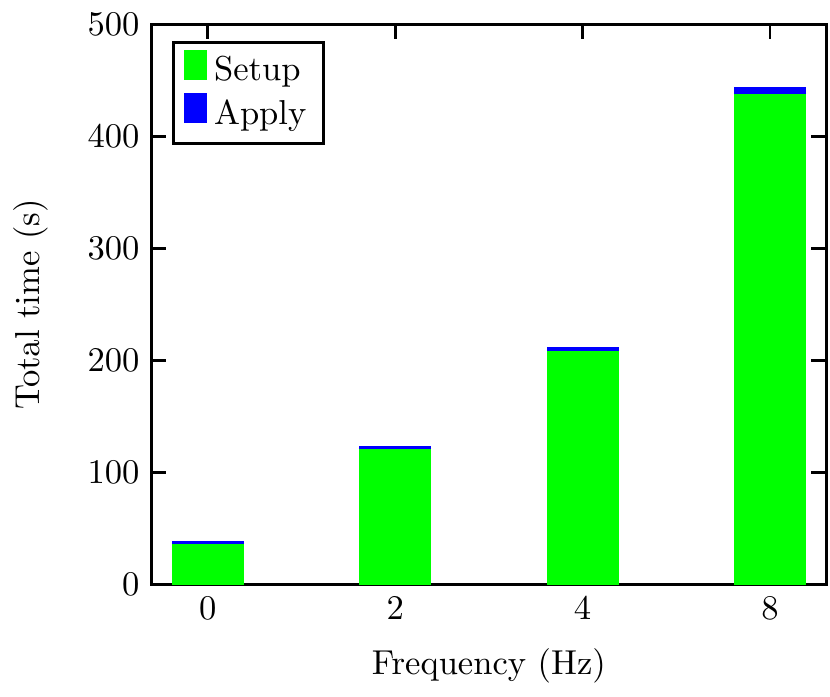}
\caption{Time requirements as a function of frequency. The high-frequency regime ($f=8$) requires the most time in both setup and application phases.}
\label{figPrec:3DHelLowHighTime}
\end{subfigure}\hspace{0.1in}
\begin{subfigure}{.5\textwidth}
\centering
\includegraphics[width=0.75\linewidth]{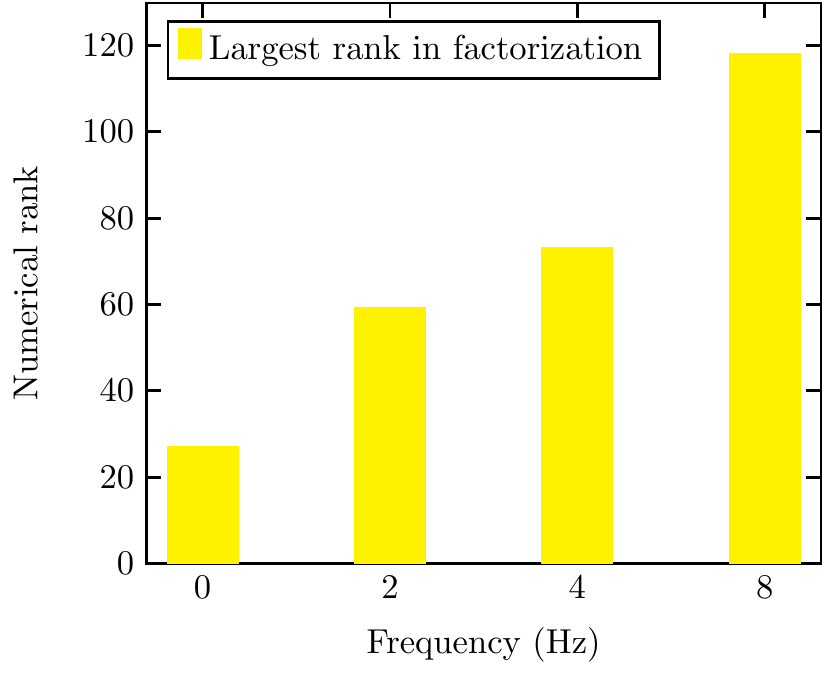}
\caption{Largest rank in factorization, the accuracy is adjusted to require less than 20 GMRES iterations. The high-frequency case ($f=8$) requires the largest numerical rank.}
\label{figPrec:3DHelLowHighRanks}
\end{subfigure}
\begin{subfigure}{0.5\textwidth}
\centering
\includegraphics[width=0.75\linewidth]{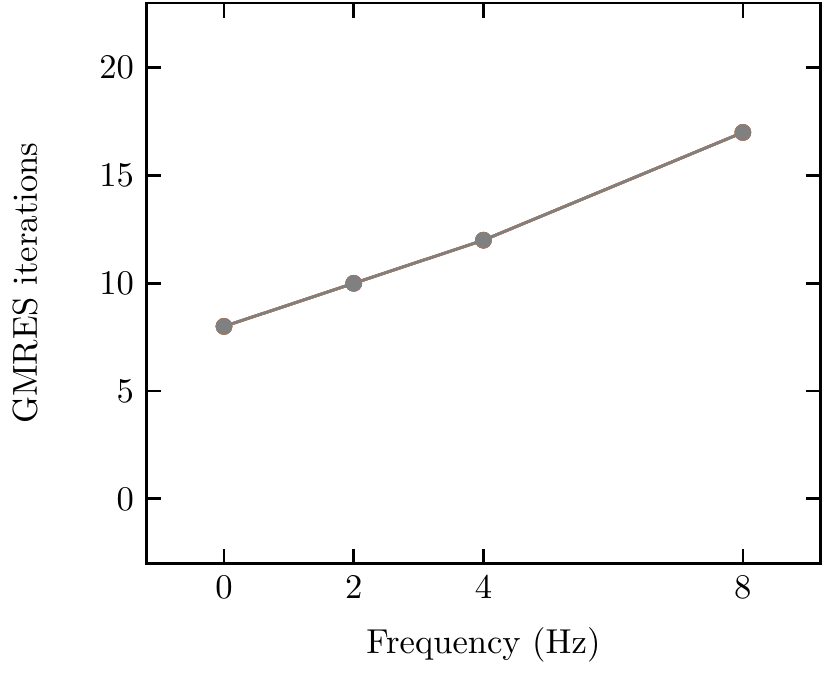}
\caption{Number of iterations as a function of frequency. The high-frequency regime ($f=8$) requires the largest number of iterations.}
\label{figPrec:3DHelLowHighIters}
\end{subfigure}\hspace{0.1in}
\begin{subfigure}{.5\textwidth}
\centering
\includegraphics[width=0.75\linewidth]{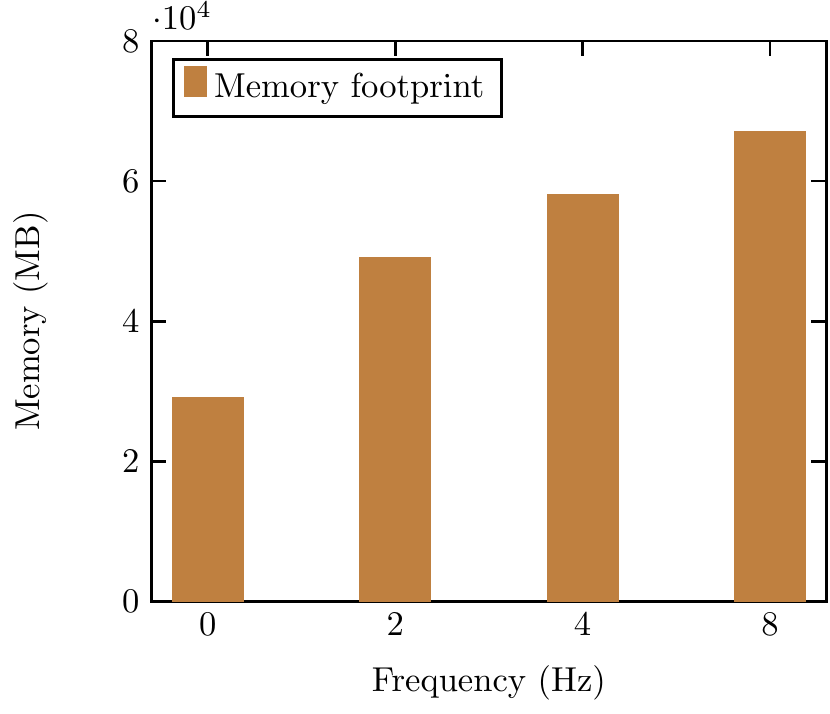}
\caption{Memory requirements as a function of frequency. The high-frequency regime ($f=8$) exhibits the largest memory footprint.}
\label{figPrec:3DHelLowHighMem}
\end{subfigure}\hspace{0.1in}
\caption{Preconditioner performance for the Helmholtz equation in a heterogeneous medium discretized with $N=128^3$ degrees of freedom at increasing frequencies. The problem with $f=0$ Hz represents a constant-coefficient Poisson problem, while $f=8$ Hz represents a moderately high-frequency Helmholtz problem.}
\label{figPrec:figPrec:3DHelLowHigh}
\end{figure}

Figure \ref{figPrec:3DHelLowHighTime} shows an apparent increase in both the setup and application phases of the preconditioner as a function of the frequency $f$. The growth in the setup is mainly due to the higher numerical ranks required to meet the upper limit of 20 iterations, as shown in Figure \ref{figPrec:3DHelLowHighRanks}. The increase in the application phase is due to both an increase in ranks, and an increase in the indefiniteness of the problem due to a higher frequency; as is evident from the growth in the number of required iterations depicted in Figure \ref{figPrec:3DHelLowHighIters}. Finally as illustrated in figure \ref{figPrec:3DHelLowHighMem}, the memory footprint also increases with the frequency as a consequence of higher ranks.

To give a comparison with traditional techniques, Table \ref{table:comparisonsHelmholtz} shows the number of iterations that GMRES without preconditioner, the incomplete LU factorization, and algebraic multigrid as preconditioner require. For this problem type, the ACR preconditioner was the only method that was able to solve all the problems under consideration.

\begin{table}[H]
\centering
\begin{tabular}{|c|c|c|c|c|c|c|c|c|}
\hline
\multirow{2}{*}{Frequency} & \multicolumn{2}{c|}{GMRES + No Prec.} & \multicolumn{2}{c|}{GMRES + ILU(30)} & \multicolumn{2}{c|}{GMRES + AMG} & \multicolumn{2}{c|}{GMRES + ACR} \\ \cline{2-9} 
& Iterations           & Solve          & Iterations          & Solve          & Iterations        & Solve        & Iterations        & Solve        \\ \hline
0                          & 582                  & 0.30           & 207                 & 201.97         & 7                 & 0.86         & 8                 & 0.93         \\ \hline
2                          & 100,000+             & -              & 100,000+            & -              & 80                & 7.22         & 10                & 1.55         \\ \hline
4                          & 100,000+             & -              & 100,000+            & -              & 100,000+          & -            & 12                & 2.13         \\ \hline
8                          & 100,000+             & -              & 100,000+            & -              & 100,000+          & -            & 16                & 3.70          \\ \hline
\end{tabular}
\caption{Number of iterations and solve time for the solution of a sequence of increasingly indefinite Helmholtz problems $N=128^3$ with a variable coefficient. Methods under consideration include the incomplete LU factorization (ILU), algebraic multigrid (AMG), and accelerated cyclic reduction (ACR).}
\label{table:comparisonsHelmholtz}
\end{table}

\subsection{Operation count and memory footprint}
\label{sec:HelmEstimates}

\begin{figure}[H]
\begin{subfigure}{.5\textwidth}
\centering
\includegraphics[width=0.75\linewidth]{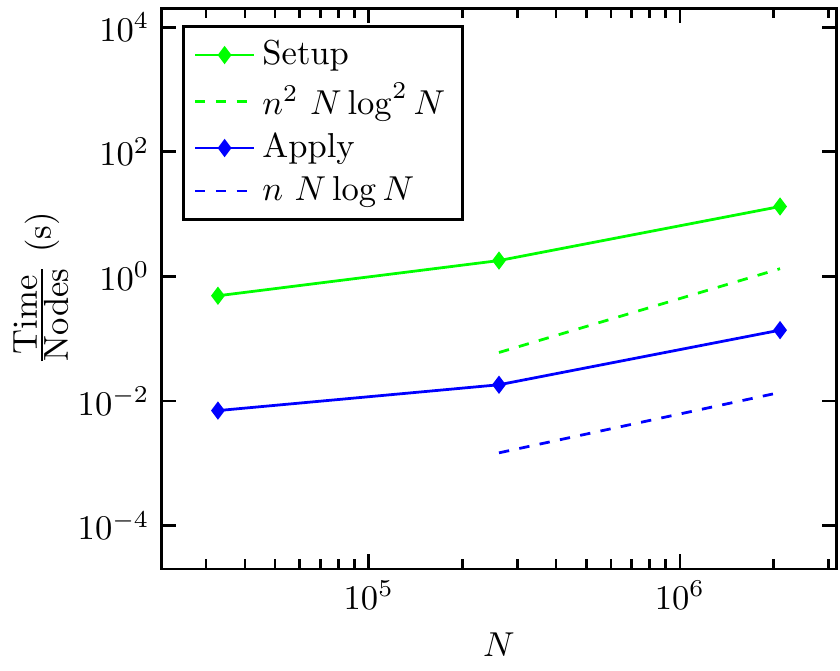}
\caption{Comparison of the preconditioner setup and application with their respective theoretical estimates.}
\label{figPrec:3DHelmComplexityTime}
\end{subfigure}\hspace{0.1in}
\begin{subfigure}{.5\textwidth}
\centering
\includegraphics[width=0.75\linewidth]{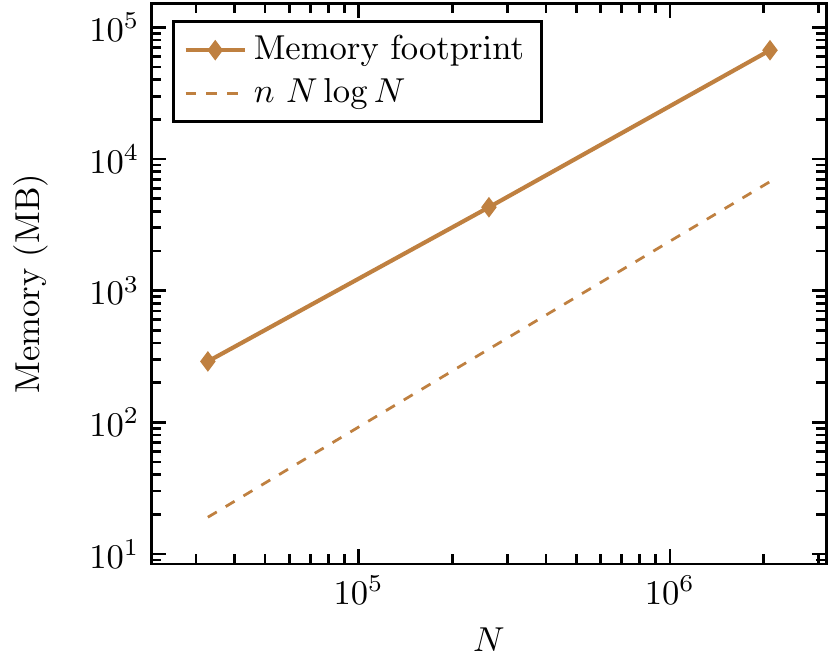}
\caption{Comparison of the preconditioner memory footprint with its theoretical estimate.}
\label{figPrec:3DHelmComplexityMemory}
\end{subfigure}
\caption{Measured performance and memory footprint for the solution of a sequence high-frequency Helmholtz problems in heterogeneous medium, discretized at 12 points per wavelength. On average, the rank of the low-rank blocks of the ACR preconditioner grows slower than $O(n)$.}
\label{figPrec:scaleNHelm}
\end{figure}

As the previous experiments show, the high-frequency Helmholtz regime is where the highest numerical ranks are required. Therefore, it is of interest to show how the computations behave asymptotically as the problem size increases considering the estimate $k \sim O(n)$ \cite{chandrasekaran2010}. Figure \ref{figPrec:3DHelmComplexityTime} shows a comparison of the preconditioner setup with the $O(n^2~N\log^2 N)$ estimate, and the preconditioner application with respect to the $O(n~N\log N)$ estimate. Figure \ref{figPrec:3DHelmComplexityMemory} shows the memory footprint of the preconditioner with respect to the estimate $O(n~N\log N)$. Table \ref{tableranks} shows a fair agreement to the Chandrasekaran et al. estimate on the largest rank growth of the factorization, however, the average rank on the low-rank blocks of the ACR preconditioner grows slower than $k \sim O(n)$, which is reflected by a slightly lower than predicted memory consumption and setup time. The ACR preconditioner does not use the HSS format or a weak admissibility condition which results in off-diagonal blocks with large rank, but rather a standard admissibility condition that allows a more refined structure of the $\mathcal{H}$-matrix blocks, as discussed in Section \ref{sec:tuningETA}, and shown in Figure \ref{fig:etaExperiment}.

\begin{table}[ht!]
\centering
\begin{tabular}{|c|c|c|}
\hline
\textbf{$N$} & \textbf{Largest rank} & \textbf{Average rank} \\ \hline
$32^3$       & 25                    & 16                    \\ \hline
$64^3$       & 59                    & 32                    \\ \hline
$128^3$      & 118                   & 36                    \\ \hline
\end{tabular}
\caption{Rank growth statistics for a sequence high-frequency Helmholtz problems in heterogeneous medium, discretized at 12 points per wavelength.}
\label{tableranks}
\end{table}

\newpage
\section{Concluding remarks}
We presented a robust and scalable preconditioner based on the cyclic reduction method and hierarchical matrices in a distributed memory environment. The preconditioner relies on a block tridiagonal structure that commonly arise from the discretization of elliptic operators with variable-coefficient.

The preconditioner setup is based on a red-black ordering for which, if a 3D grid is considered, the ordering divides the grid into planes. These planes represent block rows of the original linear system, which are represented with a hierarchical matrix, in the $\mathcal{H}$ format, and its structure is defined using a binary spatial partitioning of the planar grid sections, employing a standard admissibility criterion that controls the rank of individual low-rank blocks.

The concurrency features of ACR constitute one of its strengths. The regularity of the decomposition allows a predictable load balance. The parallel features are demonstrated via the companion implementation in a distributed memory environment with numerical experiments that study the strong and weak scalability of the method. In our current implementation, concurrency at node level involves task-based parallelism of the hierarchical matrix arithmetic operations involved in the computation of the Schur complement and its evaluation. In future work, we plan on developing a set of distributed-memory hierarchical matrix operations that can exploit a larger set of processors to accelerate the setup phase of the preconditioner.

We demonstrated over a range of problem sizes and parameters that the preconditioner can tackle a broad class of problems that lack definiteness, such as the indefinite high-frequency Helmholtz equation in heterogeneous media, or lack symmetry, such as the convection-diffusion equation with a recirculating flow.

Since the accuracy of the $\mathcal{H}$-matrix approximations and their arithmetic operations can be tuned, it was demonstrated that the preconditioner could control the number of Krylov iterations. Furthermore, we discuss how these parameters can be used to optimize performance and memory consumption via comparisons with Krylov methods with established preconditioners such as algebraic multigrid and the incomplete LU factorization, and with direct solvers that perform a complete LU factorization.

As expected from all hierarchical low-rank approximations methods, the key to performance, and memory economy, is largely based on achieving an approximation with \textit{low} rank; i.e. an efficient compression into a data-sparse format where $k$ (the rank) is much less than $n$ (the size of the block to be approximated). Numerical examples demonstrate that the required ranks agree with theoretical estimates and that for problems larger than a dozen of millions of unknowns the strong admissibility condition required less memory than the alternative (weak) admissibility condition.

\section{Acknowledgements}
We thank the editors and the reviewers for their time and comments during the review process of this work. Support from the KAUST Supercomputing Laboratory and access to Shaheen Cray XC40 is gratefully acknowledged.

\bibliographystyle{elsarticle-num}
\bibliography{jcam_acr_arxiv}

\vspace{1in}
\textbf{BibTeX entry of this article:}
\begin{verbatim}
@article{Chavez2017,
author = "Gustavo Ch{\'a}vez and George Turkiyyah and Stefano Zampini and David Keyes",
title = "Parallel accelerated cyclic reduction preconditioner for three-dimensional
elliptic \{PDEs\} with variable coefficients ",
journal = "Journal of Computational and Applied Mathematics",
year = "2017",
issn = "0377-0427",
doi = "https://doi.org/10.1016/j.cam.2017.11.035",
url = "https://www.sciencedirect.com/science/article/pii/S0377042717305952",
}
\end{verbatim}

\end{document}